\documentclass[a4paper,11pt]{article}

\usepackage[colorlinks=false,pdfborderstyle={/W 1}]{hyperref}

\def\arXiv#1#2{\href{http://arxiv.org/abs/#1}{{\tt arXiv:#1 [#2]}}} 
\def\arXivo#1{\href{http://arxiv.org/abs/#1}{{\tt [arXiv:#1]}}} 

\usepackage{import}
\usepackage{amsmath}
\usepackage{amssymb}
\usepackage{amsthm}
\usepackage{amsfonts}
\usepackage{graphicx}
\usepackage{mathrsfs}
\usepackage[margin=1in]{geometry}
\usepackage{tabularx}
  \ifx\LabelFigloaded\MYundefined\relax
  \else
   \fi

  \def\LabelFigloaded{\relax}

  \chardef\LabelFigCatAt\the\catcode`\@
  \catcode`\@=11

 \let\LabelFigwlog@ld\wlog
 \def\wlog#1{\relax}

 \ifx\\\MYundefined@
    \let\\\relax
 \fi

  \def\ms@g{\immediate\write16}

 \def\N@wif{\csname newif\endcsname }
 \def\Temp@ {\N@wif\ifIN@}
 \ifx\INN@\MYundefined@
    \else \let\Temp@\relax
 \fi
 \Temp@

  \def\IN@{\expandafter\INN@\expandafter}
  \long\def\INN@0#1@#2@{\long\def\NI@##1#1##2##3\ENDNI@
    {\ifx\m@rker##2\IN@false\else\IN@true\fi}%
     \expandafter\NI@#2@@#1\m@rker\ENDNI@}
  \def\m@rker{\m@@rker}
 
  \newtoks\Initialtoks@  \newtoks\Terminaltoks@
  \def\SPLIT@{\expandafter\SPLITT@\expandafter}
  \def\SPLITT@0#1@#2@{\def\TTILPS@##1#1##2@{%
     \Initialtoks@{##1}\Terminaltoks@{##2}}\expandafter\TTILPS@#2@}

 \def\Shifted@@#1#2#3{\setbox0=\hbox{#3}%
   \raise -\dp0\vbox {\kern-#2%
       \hbox {\kern#1\unhbox0\kern-#1}%
           \kern#2}}

 \newcount\gridcount
 \newbox\auxGridbox@ \newbox\hGridbox@ \newbox\vGridbox@
 \newbox\Labelbox@ \newbox\auxLabelbox@
 \newbox\Coordinatebox@
 \newtoks\Labeltoks@
 \newdimen\Wdd@ \newdimen\Htt@
 \newdimen\Wddd@ \newdimen\Httt@
 
 \def\Wr@{\immediate\write16}

 \newdimen\GL@wd
 \def\GridLineWidth#1{\GL@wd=#1}

 \def\gobble#1{}
 \def\EdgeErr@{\Wr@{}%
      \Wr@{\string\Edges\space argument
      1, 10, 100 or 1000 please\string!}%
      }

 \newcount\Edgect@

 \def\Sweepup#1\endSweepup{}

 \def\SetEdges@{%
    \edef\Zr@@s{\expandafter\gobble\number\Edgect@\empty}%
        \count255=0\Zr@@s\relax
        \ifnum\count255=\z@\else\EdgeErr@\show\tailtest\fi
        \count255=1\Zr@@s\relax
        \ifnum\count255=\Edgect@\relax\else\EdgeErr@\show\leadtest\fi
    \EdgGl@b\edef\Zr@s{\expandafter\gobble\Zr@@s\empty}
    \ifnum\Edgect@>\@ne\relax\EdgGl@b\let\L@Dc\empty
        \else\EdgGl@b\edef\L@Dc{\string.}\fi
    \ifnum\Edgect@>\@ne\relax
        \EdgGl@b\edef\Edgescale@##1{\divide##1 by \Edgect@}%
        \else\EdgGl@b\edef\Edgescale@##1{}\fi
    }

 \def\Edges#1{\Edgect@=#1\relax
     \let\EdgGl@b\global \SetEdges@}

 \Edges{1}

 \def\hhrule{\hrule height \GL@wd\vskip-.\GL@wd}

 \def\hRule@{%
   \advance\gridcount -2%
   \vfil\hhrule\vfil
   \llap{\smash{\raise -2.5pt
     \hbox{\L@Dc\number\gridcount\Zr@s\kern2pt}}}%
   \hhrule
   }

\def\vvrule{\vrule width \GL@wd \kern-\GL@wd}

 \def\vRule@{\advance\gridcount 2%
   \hfil\vvrule\hfil
   \setbox\auxGridbox@=\vbox to 0pt
      {\vskip \Htt@\vskip 2pt
        \hbox to 0pt{\hss\L@Dc\number\gridcount\Zr@s\hss}\vss}%
      \wd\auxGridbox@=0pt \box\auxGridbox@
   \vvrule
   }

 \def\PlaceGrid@@{\gridcount=10 
  \setbox\hGridbox@=\hbox{%
        \hbox{%
             \hskip-.4pt\vrule
             \vbox to \Htt@{%
               \offinterlineskip\parindent=\z@\relax
               \hbox to \Wdd@{\hfil}
               \hRule@\hRule@\hRule@\hRule@
               \vfil\hhrule\vfil}%
             \vrule\hskip-.4pt}
    }%
  \gridcount=0%
  \setbox\vGridbox@=\hbox{%
      \vbox{\offinterlineskip\parindent=0pt\hsize=0pt
         \vskip-.4pt\hrule%
         \hbox to \Wdd@{%
                 \vtop to \Htt@{\vfil}%
                 \vRule@\vRule@\vRule@\vRule@
                 \hfil\vvrule\hfil}%
         \hrule\vskip-.4pt}}%
  \wd\hGridbox@=0pt\ht\hGridbox@=0pt
  \wd\vGridbox@=0pt\ht\vGridbox@=0pt
  \hbox{\box\hGridbox@\box\vGridbox@}%
  }

 \def\LabelsGlobal{\def\LabGl@b{\global}}
 \def\LabelsLocal{\def\LabGl@b{}}
 \LabelsGlobal 

 \def\SetLabels#1\endSetLabels{%
   \LabGl@b\Labeltoks@={#1()\\}%
   }

 \LabGl@b\Labeltoks@={()\\}

 \def\ShowGrid{\LabGl@b\let\PlaceGrid@\PlaceGrid@@}
 \def\HideGrid{\LabGl@b\let\PlaceGrid@\relax}
 \def\Grids{\ShowGrid\LabGl@b\let\GridSwitch@\ShowGrid}
 \def\noGrids{\HideGrid\LabGl@b\let\GridSwitch@\HideGrid}

 \noGrids

 \def\bAdjust@@{%
     \setbox\auxLabelbox@=\hbox{\raise \dp\auxLabelbox@
            \box\auxLabelbox@}}
 \def\bAdjust@{\let\vAdjust@\bAdjust@@}

 \def\eAdjust@@{\dimen0=-.5\ht\auxLabelbox@
     \advance\dimen0 by .5\dp\auxLabelbox@
     \setbox\auxLabelbox@=
            \hbox{\raise\dimen0\box\auxLabelbox@}}
 \def\eAdjust@{\let\vAdjust@\eAdjust@@}

 \def\tAdjust@@{%
     \setbox\auxLabelbox@=\hbox{\raise-\ht\auxLabelbox@
            \box\auxLabelbox@}}
 \def\tAdjust@{\let\vAdjust@\tAdjust@@}

 \let\vAdjust@\relax

 \def\lAdjust@{\let\hAdjust@\rlap}
 \def\rAdjust@{\let\hAdjust@\llap}

 \let\hAdjust@\relax\let\vAdjust@\relax

 \def\FetchLabel@#1(#2)#3\\{%
     \IN@0#2@@\ifIN@
        \setbox0=\hbox{\ignorespaces#1#3\unskip}%
        \ifdim\wd0>0pt
           \ms@g{}%
           \ms@g{ !!! Bad label(s)? !!!}%
           \message{ #1(#2)#3}%
        \fi
        \def\LabelMole@##1\endFetchLabel@{%
            \IN@0()\\@##1@%
            \ifIN@\def\Temp@{\FetchLabel@##1\endFetchLabel@}%
            \else\def\Temp@{}%
            \fi
            \Temp@
           }%
     \else
       \ignorespaces#1\unskip
       \setbox\auxLabelbox@=%
         \hbox to 0pt{\hss\ignorespaces\hAdjust@
          {\ignorespaces#3\unskip}\hss}%
       \vAdjust@
       \let\hAdjust@\relax\let\vAdjust@\relax
       \AugmentLabelBox@@{#2}%
       \ht\Labelbox@=0pt\dp\Labelbox@=0pt
       \let\LabelMole@\FetchLabel@%
     \fi\LabelMole@}

 \newtoks\XYSep@ 
 \def\SetXYSeparator#1{%
     \IN@0#1@@\ifIN@\XYSep@{*}%
     \else
     \XYSep@{#1}%
     \fi
     }

 \SetXYSeparator*

 \def\AugmentLabelBox@@#1{%
     \IN@0\the\XYSep@ @#1@\ifIN@
       \SPLIT@0\the\XYSep@ @#1@%
       \setbox\Labelbox@=\hbox to 0pt{%
         \unhbox\Labelbox@
         \Shifted@@{\the\Initialtoks@\Wddd@}%
         {\the\Terminaltoks@\Httt@}%
         {\box\auxLabelbox@}}%
     \else
         \ms@g{}%
         \ms@g{ !!! Bad insertion point. !!!}%
         \message{ (#1\ this point was rejected.)}%
     \fi
    }

 \def\FetchOption@#1[#2]#3\endFetchOption@{%
    \def\temp{#1}
    \ifx\temp\empty
       \Edgect@=#2\relax
       \let\EdgGl@b\relax
       \SetEdges@
       \Cleaner@#3%
    \fi}

 \def\Cleaner@#1[@]{\Labeltoks@{#1}}
     
 \def\PlaceLabels@@{\mathsurround=0pt
     \def\Cr@{\\}%
     \let\L\lAdjust@\let\R\rAdjust@
     \let\B\bAdjust@\let\E\eAdjust@\let\T\tAdjust@
     \expandafter\FetchOption@\the\Labeltoks@[@]\endFetchOption@
     \Wddd@=\Wdd@ \Edgescale@\Wddd@ 
     \Httt@=\Htt@ \Edgescale@\Httt@
     \expandafter\FetchLabel@\the\Labeltoks@\endFetchLabel@
     \box\Labelbox@
     }%

 \let \PlaceLabels@\PlaceLabels@@

 \def\AffixLabels#1{\setbox\Coordinatebox@=\hbox{#1}%
      \Wdd@=\wd\Coordinatebox@ \Htt@=\ht\Coordinatebox@
      \advance\Htt@ \dp\Coordinatebox@
      \hbox{\copy\Coordinatebox@\kern-\Wdd@ 
           \Shifted@@{0pt}{-\dp\Coordinatebox@}%
           {\PlaceLabels@\PlaceGrid@}%
           \kern\Wdd@}%
      \GridSwitch@ 
      \LabGl@b\Labeltoks@{()\\}%
      }
 
   \let\wlog\LabelFigwlog@ld   
   \catcode`\@=\LabelFigCatAt  

                                By

              Raymond S\'eroul <A18645@FRCCSC21.BITNET>
                                and 
              Laurent Siebenmann <lcs@topo.math.u-psud.fr>
    
              VERSIONS: July 1991, Oct 1991, Jan 1992, July 1992

INTRODUCTION

      This labelling package is intended for TeX users who
rely on non-TeX sources for for their graphics inserts.  It
provides means for adding TeX labels to such inserts with a
minimum of fuss. 

       For most labels, TeX users have in the past found it
reasonably convenient to rely on non-TeX sources. Typical
occasions when an inescapable need for TeX labels seemed to
arise are

 (a) when the graphics program lacks certain exotic or complex
mathematical symbols

 (b) when the very highest typographical quality is wanted for the
labels

 (c) when labels included with the graphics fail to print, 
 and you cannot figure out why (cf. boxedeps.doc).  The labels
 provided by labelfig.tex are 100

       Since this package first appeared, many users, who in the
past scarcely dreamed of using TeX labels, have come to use
nothing but.  So it is now appropriate to add

Intoxication Warning:  TeX labels may be addictive and expensive. 

     If you have a fast preview you may disagree, and even find
that this package provides an agreeable paste-up environment; see
extra applications at end.

     Note to publishers: It is possible and convenient to ultimately
export the TeX labels produced by labelfig.tex to become an integral
part of the EPS file. This is often desired by a publisher who typically
uses an "upmarket" graphics or page layout program, with which the
staff is skilled in perfecting figures.  See Appendix I for
a recipe.

     The authors are grateful to Patrick Ion of Math Reviews for
helpful comments and encouragement.

BASIC INSTRUCTIONS

    After reading in the macro file using

preview or proof your figure with a coordinate grid printed on
top, by typing the following:

    \ShowGrid  
    \AffixLabels{<the graphics insertion>}

Here <the graphics insertion> is what you would type to insert
the graphics object alone without the grid.  This must provide
for the space around it. For example <the graphics insertion>
might well be \BoxedEPSF{MyFigure scaled 700} using the
boxedeps.tex macro package (from same source); this provides a
TeX box containing the encapsulated PostScript insert specified by
the file MyFigure. \AffixLabels{...} provides the grid (supposing
\ShowGrid is present) and later, once you have specified labels
using the grid, it will "tack on" the labels.

     The grid is a sort of (usually elongated) checkerboard of
ten rows and ten columns and its (internal) partitions are by
default numbered  .1, ... ,.9  both horizontally (X-coordinate
running left to right) and vertically (Y-coordinate running bottom
to top).  Thus the points enclosed by the grid correspond to the
points of the unit square in the cartesian "X-Y" plane, the lower
left corner corresponding to the origin (0,0).  By extrapolation,
the full page corresponds to a larger rectangle in the plane.

     These coordinates serve to position labels as follows.
Before the \AffixLabels{...} command type label specifications:

  \SetLabels
   (<X-coordinate>*<Y-coordinate>) <first label> \\
   .
   .
   .
   (<X-coordinate>*<Y-coordinate>)  <last label> \\
  \endSetLabels

Each row specifies one label and is terminated by \\.  In each
row, the position indicator comes first; it is written as a
standard cartesian point except that the X- and Y- coordinates
are separated by * rather than a comma because TeX allows a
comma as decimal point. There are no dimension units to specify
as the unit is the grid itself.

     By default, this cartesian point specifies where the middle
of the baseline of the label will be located.  However if you precede
the point by \L [or \R] the left [or right] edge of the baseline will
be located there. Similarly you may also precede the point by \T, \E,
or \B to vertically align the top equator or bottom of the label box
at the specified point.  This gives nine standard positions of
the label with respect to the insertion point --- corresponding to
the eight principle points of the compas and the center

                     \L\T     \T      \R\T

                     \L\E     \E      \R\E

                     \L\B     \B      \R\B

But this neglects the default "baseline" level of TeX,
giving potentially three more positions

                     \L    <no tag>   \R

For text, the baseline level is often the preferred. Its relation to
the others is variable. It will often coincide with the bottom level,
as happens for "X".  But it is often distinct, as for "g", in which
case you have in all 12 distinct positions rather than 9.

     It is convenient to think of this specification of label
position as attaching the label by a thumb-tack to the coordinate
grid. There are up to twelve positions of the thumb-tack on the
label, while the position of the thumb-tack on the coordinate grid is
arbitrary.  Normally, one choses the position of the thumb-tack on
the label to be the one that is the closest to the item being
labeled.  There are good reasons for this "rule of thumb":

   (a)  It facilitates correct positioning at first try.

   (b)  If the scale of the figure must be altered after labels
have been affixed, the labels have a good chance of remaining well
positioned.

   (c)  The visible grid need not extend beyond the "bounding box"
for the figure, because the best preferred position is always
(at least almost) within the bounding box .

The second reason is particularly important. Indeed it often
happens that scale has to be altered after labelling begins, in
order to either provide space for the labels, or to adjust
proportions between the labels and the figure.  (The size of labels
is unaffected by scaling.)

     Here is an artificial but self-contained test which uses
TeX rules to make a graphics object.

TEST

    Do not skip this!

 \def\FrameIt#1{\hbox{\vrule$\vcenter {\hrule\kern3pt%
             \hbox {\kern3pt #1\kern3pt}%
               \kern3pt\hrule}$\relax\vrule}}

 \def\Caption#1#2{\FrameIt{%
       \vtop {\hsize=#1\relax \parindent=0pt
         \leftskip=0pt \rightskip=0pt plus15pt
         \parfillskip=0pt
         \lineskip=1pt\baselineskip=0pt
         #2}}}

 \def\FirstQuadrant{\hbox to 100pt{\vrule\vbox to 100pt{%
        \hbox to 100pt{\hfil}\vfil\hrule}\hss}}

  \SetLabels
    \R(.5*.2) $\zeta\,\cdot$\\
    (.9*-.10) $\xi$\\
    \R(-.03*.9) $\eta$\\
    \T(.5*.9) \Caption{70pt}{%
          \it The norm of
          $g(\xi+i\eta)$ is indicated on
          contours of this invisible surface.}\\
  \endSetLabels

  \AffixLabels{\FirstQuadrant} \end

  Note that the coordinates to use for labels are indicated on the
edges of the grid (when visible) corresponding to the conventional
x- and y- axes of the Cartesian plane. By default the grid is
1-by-1. However, by the command \Edges{100}, you can change this
to 100-by-100 and many users find this alternative most
convenient. Place the command \Edges{...} in your style file (or
header) since its effect is is global. Other possible edge values
are 10 and 1000.

  If you use the command \Edges{...} at all, do so with care.  For
if you accidentally delete an \Edges{...} command your labels will
abruptly be badly misplaced and may logically but mysteriously
generate "dimension too big" errors under TeX and "off page" errors
under your driver.  

  You can dictate the edgescale for an individual figure by giving
the scale in brackets immediately after \SetLabels.  Thus, to
import into an article using say \Edge{100} a figure labelled using
another edgescale, say the original 1-by-1 default, you can use
\SetLabels[1]...\endSetLabels.

GETTING IT DOWN PAT

     Complicated labeling deserves the same respect as
complicated mathematics.  Do not expect it to come out perfect the
first time!  What is needed in either case is a mechanism to
repeatedly typeset troublesome pieces.

     One mechanism is always available.  One does complicated
labelling in a separate "test" file involving just the figure being
labelled;  a texpert will know how to \dump TeX's current state as
a temporary format that restarts rapidly at each retry.  Usually,
one then pastes the completed labelled figure back into the main
TeX file, but, of course, one can also \input it as an auxiliary
file.

     If you do not have a TeXpert at handy, here is a first
approximation to an efficient setup. By deletions reduce a copy
of your article to just a few lines before and after the figure.
Now label the figure, and finally, copy and paste the labelled
figure to the original article. Then copy the next figure to label
into this testbed and repeat. The TeXpert can improve the  speed
at which TeX starts up, by compiling a format specifically for
your article; just one caution: best NOT include in the format
ephemeral details of setup like \Set<mydriver>ArtSpecials (from
boxedeps.tex because this reads  figure dimensions which you may
change during your work session.

     An improved mechanism to repeatedly typeset troublesome
pieces is now available on the Macintosh; it is called LinoTeX;
see the same ftp sources.  It could be set up on many types
of computer.

     Before using labelfig.tex to attach labels to a graphics
object inserted using boxedeps.tex or BoxedArt.tex, make it a
firm rule to carefully adjust the bounding box using the trimming
commands of these packages, and also at least tentatively scale
and position the object. Beware of changing the grid inadvertently
after the labels have been positioned.  For example, correcting
the bounding box of a PostScript graphics object can foul up the
labels by changing the coordinate grid to which the labels are
attached. This is particularly true for the trimming  commands of
boxedeps.tex and BoxedArt.tex. However, as noted already, change
of scale is much less disruptive, and modest adjustments should be
well tolerated.

     Sometimes the labels protrude so far from the bounding box
of a figure that the figure has to be repositioned.  Best do this
by ad hoc spacing, say using \hglue and \vglue; altering the
bounding box would create a vicious circle.

     Remember that you are responsible for preventing labels
from overlapping. You are responsible for all label typography
including size and style. A label is really just about anything
that can be put in a TeX box. Note that spaces at the beginning
and end of labels will normally be suppressed; if you really want
them you must protect them with TeX braces.

     This package temporarily sets the \mathsurround parameter
of TeX to zero  while the labels are being affixed. This is done
because nonzero \mathsurround space would influence the position
of left and right aligned labels; then, when a texpert or printer
modifies mathsurround, diagram labeling might be disastrously
altered. There is a small price to pay involving labels that are
formatted as caption boxes including mathematics: you  may want or
need to specify an explicit mathsurround space within the caption
box; it will not influence anything outside.

     Those hostile to the use of * as separator between
the X and Y coordinates of label insertion points, are free to
impose another using \SetXYSeparator{<the new separator>}.  
Americans may prefer "," to "*" since they never use a 
comma as a decimal point; on the other hand, * may be more visible.

APPENDIX (I)  MERGING labelfig.tex LABELS INTO AN EPSF GRAPHICS OBJECT.

     As promised in the introduction, here is a recipe useful for
publishers. It works at least on Macintosh and at least for vectorized
graphics and Adobe type1 fonts.  (There is surely a similar recipe for
PCs under MSWindows.)

 (a)  Use boxedeps.tex utility to integrate the figure given by the eps
file, "x.eps" say, with a visible frame around it.  See
\ShowDisplacementBoxes command in boxedeps.tex.  To get precise results
automatically it is important to use the \Trim... commands of
boxedeps.tex making the "DisplacementBox" neatly fit the figure.

 (b)  Use the TeX printer driver and LaserWriter (versions >= 8.1.1) to
export to an EPSF the DVI page containing the integrated, labelled
figure. You now have an EPS file  "xx.eps"  that contains too much, and at
the wrong scale, and at wrong position.

 (c)  Convert the EPSF to an Adode Illustrator format EPSF using
the shareware utility called epsConvert by Sam Weiss
1993-- (currently $25).

 (d)  In Illustrator (or a compatible program), group the labels and the
"DisplacementBox"; copy them to the clipboard and paste them into "x.ps".
This step requires that all the label fonts be "visible to the Macintosh.

 (e)  Translate and scale the pasted group consisting of the labels plus
the "DisplacementBox" so as to make the "DisplacementBox" the bounding
box of (labelless) figure represented by "x.eps".  At this point the
labels will be correctly placed on the figure "x.eps".

 (f)  Ungroup and delete the "DisplacementBox".  The result is the
desired single EPS file, "x+.eps" say, It contains the original figure
plus its labels.  

     Using grouping and ungrouping appropriately in "x+.eps", a
publisher's staff can very efficiently improve label positions etc.

APPENDIX II)  SOME EXOTIC APPLICATIONS

     The grid of labelfig.tex is analogous to a light-table in
classical page makeup with wax or latex glue.  In principle, you
can use it to compose any page from its indivisible parts.  This
even has some of the artisanal charm of classical paste-up
provided you have a fast screen preview to make the process
"interactive".

     In practice labelfig.tex is a tool for nonstandard jobs.
Here are a few going beyond the labelling already discussed.

(I)  GRAPHICS INTEGRATION.

     This is accomplished by treating the imported graphics
objects as labels.  The underlying graphics object is then
typically an empty  \vbox to <dimension>{\vfill} in a TeX
\midinsert...\endinsert construction.  A label line
might be of the form

   (.1*.1) \\

The exact form of the special command varies from driver to
driver.  However, in the case of encapsulated PostScript graphics
(EPSF norm), by relying on boxedeps.tex, one can have the
following standard syntax (independant of driver  (see
boxedeps.doc for details.
  
  (.1*.1) \BoxedEPSF{MyFigure scaled <scale in mils>}\\

This may be slow since it requires TeX to read the PostScript
file to read bounding box using many complex macros.  So you
may want to try

  (.1*.1) \EPSFSpecial{MyFigure}{<scale in mils>}\\

which is fast and driver independant, but it squashes the
bounding box, normally to its lower left corner.

     Similarly for graphics of the Macintosh PICT norm ---
using BoxedArt.tex (same sources) in place of boxedeps.tex.

     This approach to integration is to be recommended when
one is assembling a composite graphics object.

 (II)  COMMUTATIVE DIAGRAM ENHANCEMENT

     Commutative diagrams or arrays of mathematical objects
connected by arrows of various sorts are common in mathematics.
The mathematical objects require the use of TeX.  Recently TeX
acquired a good collection of arrows of all slopes --- that of
LamSTeX --- plus pwerful macros to build the diagrams.

     However, even the LamSTeX collection is often
inadequate; it lacks for example double shafted arrows, dotted
arrows and curved arrows. Fortunately it is possible to produce
such arrows on an individual basis using sophisticated graphics
programs such as Illustrator and AldusFreehand (both serving
the EPSF norm) or using Metafont (with its public domain norm).
Since the creation of each new arrow is a work of love, you
probably want to limit the number of arrows by using LamSTeX
for most arrows. The 40K commutative diagram module of LamSTeX
has been adapted to work with AmSTeX and a copy may be posted
with LabelFig and related files. Unfortunately no one has yet
offered a version that works with Plain TeX or LaTeX.

       Suffice it here to say that when the exotic arrow has
been somehow imported into TeX, labelfig.tex treats it as a
label that one affixes to the commutative diagram.  Two other
steps will be treated in separate notes, namely the matter of
extracting the dimension specifications for the arrow and the
construction of the arrow --- for these steps are far from
unique and often depend intimately on your computer environment. 
Notes for the Macintosh-Textures-Illustrator combination are
found in the file ExoticArrows.doc.

 (III) NESTING 

Ingenuity pays off in exploiting labelfig.tex. One can
mix graphics and typography quite freely.  labelfig.tex is good
for freeform or overlapping arrangements, while boxedeps.tex (or
BoxedArt.tex) is best for regimented non-overlapping
arrangements --- and the two can be combined.

     The default behavior of labelfig.tex is not ideal 
for nesting objects, because to prevent trouble for beginners
the register for labels is globally cleared when \AffixLabels
concludes.  But there are switches available

      \LabelsGlobal      \LabelsLocal

which change this.  To understand this, extend the above test 
by something like:

 \LabelsLocal

 \SetLabels
    (.5*.5) AAA\\
 \endSetLabels

 {
 \SetLabels
    (.5*.5) ZZZ\\
 \endSetLabels
   \AffixLabels{\FirstQuadrant}
 }

   \AffixLabels{\FirstQuadrant}

     There are however potential pitfalls.  Neither
labelfig.tex nor boxedeps.tex has been tested under extreme
conditions. Problems may occur if their procedures are
indiscriminately nested. For boxedeps.tex (not labelfig.tex)
there is a precise cause for worry, namely many of its
variables are "global", which means that TeX braces will not
provide the protection one might expect.

COMMAND SUMMARY FOR labelfig.tex

  Here [...] means optional (one or zero)
       [...]* means any number of such constructs

  \SetLabels
    [[<P>](<X><Sep><Y>) <label> \\]*
  \endSetLabels
  \ShowGrid  
  \AffixLabels{<the figure>}

   --- <P> is tack position, one of eleven or empty
              order irrelevant

                   \L\T      \T      \R\T

                   \L\E      \E      \R\E

                     \L               \R

                   \L\B      \B      \R\B

   --- (<X><Sep><Y>) insertion point;
  <Sep> is separator, = * by default;
  \SetXYSeparator{<Sep>} changes it.
   <X> and <Y> are real numbers

  --- <label> a label to attach 

  --- <the figure> the figure to label 

  \GlobalLabels (default)     
  \LocalLabels  setting for nested constructs.

 \Grids makes ALL grids appear; \HideGrid then makes just next disappear.
 \noGrids returns to default.  The commands are always global.

 \GridLineWidth{<dimension>} adjusts width of grid lines. Default is very
small, to give "hairline" effect. If your grid lines are missing try
setting \GridLineWidth{1pt}.

 \Edges#1 globally changes the edge size of all grids to the numerical 
value #1, which must be 1, 10, 100, or 1000.  The default is 1.

VERSION HISTORY.
 --- Jan 1993: \Edges#1 and [??] option after \SetLabels
 --- July 1992: \Grids, \noGrids, \HideGrid;
       Gridlines become hairlines; \GridLineWidth{<dimension>}.
 --- Oct 1991, Jan 1992: \SetXYSeparator{<Sep>},  \LabelsGlobal,
       \LabelsLocal.
 --- July 1991: first release

Address for bugs and other feedback:

        Raymond S\'eroul
        IREM and Lab. de Typographie Informatise
        Univ. Rene Descartes
        Strasbourg

    Tel 33-88-41-63-45
    Email:  A18645@FRCCSC21.BITNET

        Laurent Siebenmann
        Mathematique, Bat. 425,
        Univ de Paris-Sud,
        91405-Orsay,
        France

    Tel 33-1-6941-7949; 
    Email: lcs@topo.math.u-psud.fr

\newtheorem{thm}{Theorem}[section]
\newtheorem{quest}[thm]{Question}
\newtheorem{conj}[thm]{Conjecture}
\newtheorem{lemma}[thm]{Lemma}
\newtheorem{proposition}[thm]{Proposition}
\newtheorem{corollary}[thm]{Corollary}
\theoremstyle{remark}\newtheorem{rem}[thm]{Remark}
\theoremstyle{definition}
\newtheorem{definition}{Definition}[section]

\newcommand{\R}{\mathbb{R}}
\newcommand{\p}{\mathbb{P}}
\newcommand{\E}{\mathbb{E}}
\newcommand{\Z}{\mathbb{Z}}
\newcommand{\N}{\mathbb{N}}
\newcommand{\Cov}{\mathrm{Cov}}
\newcommand{\Corr}{\mathrm{Corr}}
\newcommand{\Var}{\mathrm{Var}}

\newcommand{\B}{\mathcal{B}}

\newcommand{\cU}{\mathcal{U}}
\newcommand{\cB}{\mathcal{B}}

\newcommand{\cF}{\mathcal{F}}
\newcommand{\cW}{\mathcal{W}}

\newcommand{\Piv}{\mathsf{Piv}}

\newcommand{\floor}[1]{\left\lfloor {#1} \right\rfloor}
\def\1{1\!\! 1}

\def\clue{\mathsf{clue}}
\def\sig{\mathsf{sig}}

\def\Corr{\mathrm{Corr}}
\def\Spec{\mathscr{S}}
\def\Piv{\mathscr{P}}

\def\eps{\epsilon}
\def\eps{\epsilon}
\def \Cr{\mathsf{LR}}
\def \Stab{\mathbf{Stab}}

\def \Tr{\mathsf{Tribes}}
\def \Maj{\mathsf{Maj}}
\def \Dict{\mathsf{Dict}}

\def\given{\;|\;}

\def\lora{\longrightarrow}

\numberwithin{equation}{section}
\numberwithin{figure}{section}
\def\bl{\begin{lemma}}
\def\el{\end{lemma}}
\def\bth{\begin{thm}}
\def\eth{\end{thm}}
\def\bc{\begin{corollary}}
\def\ec{\end{corollary}}
\def\bcj{\begin{conj}}
\def\ecj{\end{conj}}
\def\bpr{\begin{proposition}}
\def\epr{\end{proposition}}
\def\bde{\begin{definition}}
\def\ede{\end{definition}}
\def\beq{\begin{equation}}
\def\eeq{\end{equation}}
\def\bpf{\begin{proof}}
\def\epf{\end{proof}}
\newcommand{\comm}[1]{}

\begin{document}

\title{ Sparse reconstruction in spin systems I: iid spins}
\author{P\'al Galicza \and G\'abor Pete
}

\date{}
\maketitle
\begin{abstract}
For a sequence of Boolean functions $f_n : \{-1,1\}^{V_n} \longrightarrow \{-1,1\}$, defined on increasing configuration spaces of random inputs, we say that there is sparse reconstruction if there is a sequence of subsets $U_n \subseteq V_n$ of the coordinates satisfying $|U_n| = o(|V_n|)$  such that knowing the coordinates in $U_n$ gives us a non-vanishing amount of information about the value of $f_n$.

We first show that, if the underlying measure is a product measure, then no sparse reconstruction is possible for any sequence of transitive functions. We discuss the question in different frameworks, measuring information content in $L^2$ and with entropy. We also highlight some interesting connections with cooperative game theory. Beyond transitive functions, we show that the left-right crossing event for critical planar percolation on the square lattice does not admit sparse reconstruction either. Some of these results answer questions posed by Itai Benjamini.
\end{abstract}

\tableofcontents

\vfill\eject

\section{Introduction and main results}\label{s.intro}

Consider some random variables $X_V:=\{ X_v : v\in V \}$ on some probability space $(\Omega^V,\,\p)$, where $\p$ is not necessarily a product measure, and a function $f: \Omega^V \lora \R$, often the indicator function of an event. How much information does a subset $X_U$ of the input variables has about the output $f(X_V)$? There are several possible approaches to formulate this question precisely; here we are focusing on the following one (and the connections to some others will be reviewed in Subsection~\ref{siginf}). Depending on the measure $\p$ and the function $f$, when is it possible that knowing $X_U$ for a small but carefully chosen subset $U$ (specified in advance, independently of the values of the variables) will give enough information to estimate $f(X_V)$?

First of all, we need to measure the amount of information we gain about $f(X_V)$ by learning a subset of the coordinate values of a configuration. For a subset $U \subseteq V$, let $\mathcal{F}_U$ denote the $\sigma$-algebra generated by $X_U$.

\bde[$L^2$-clue] \label{clue}
Let $f: (\Omega^V,\,\p) \lora \R$ and $U \subseteq V$. Then,
\beq
 \clue(f\;|\;U) :=  \frac{\Var(\E[f \,|\, \mathcal{F}_U])}{\Var (f)}.
\eeq
\ede

This notion of $\clue(f\;|\;U)$ quantifies the proportion of the total variance of $f$ attributed to the variance of the function projected onto $\mathcal{F}_U$. The clue is always a number between $0$ and $1$, as an orthogonal projection can only decrease the variance. The clue satisfies a trivial but important monotonicity property: whenever $W \subseteq U$, we trivially have $\mathcal{F}_W \subseteq \mathcal{F}_U$, and thus $\clue_{f}(W)\leq \clue(f\;|\;U)$. This is again immediate from the fact that the conditional expectation is an orthogonal projection. It is also worth noting that 
\beq\label{corr2clue}
 \clue(f\;|\;U) = \frac{\Cov^2(f,\E[f \,|\, \mathcal{F}_U])}{\Var(f)\Var(\E[f \,|\, \mathcal{F}_U])} = \Corr^2(f, \E[f \,|\, \mathcal{F}_U]),
\eeq
using that $\Cov(f,\E[f \,|\, \mathcal{F}_U]) =\Var(\E[f \,|\, \mathcal{F}_U])$.

This concept first appeared under this name in \cite{GPS}, where $f$ was always the indicator function of some crossing event in critical planar percolation, and, among many other results, the following \cite[Conjecture 5.1]{BKS} was proved. Let $\p$ be the product measure in which every edge in the box $[n]^2$ of the square lattice $\Z^2$ is deleted with probability $1/2$ independently; let $f_n$ be the indicator of having a left-to-right crossing in $[n]^2$, and let $U_n$ be the set of vertical edges in $[n]^2$. Then, $\clue(f_n \given U_n) \to 0$, as $n\to\infty$. The present project was started by Itai Benjamini asking the question (personal communication): in critical planar percolation, does $\clue(f_n \given U_n) \to 0$ hold for every sequence of subsets with $|U_n| = o(n^2)$? The results of \cite{GPS} give an affirmative answer for many such sequences $\{U_n\}$, but not for all. And how about other natural Boolean functions in place of percolation events, still with iid measures $\p_n$?

\bde[Sparse reconstruction]
Consider a sequence $f_n: (\Omega^{V_n},\,\p_n) \longrightarrow \R$. We say that there is sparse reconstruction for $f_n$ w.r.t.~$\p_n$ if there is a sequence of subsets  $U_n \subseteq V_n$  with $\lim_{n} \frac{\left|U_n\right|}{\left|V_n\right|}=0 $ such that
$
 \liminf_{n}{ \clue(f_n\;|\;U_n)}>0\,.
$
\ede

In this paper, our main focus will be on product measures $\p_n$. If $f_n$ depends only on a small proportion of the variables (e.g., dictators and juntas), then sparse reconstruction is obviously possible. But what happens if all the variables play an equal role, i.e., if there is some transitive group action $\Gamma_n \curvearrowright V_n$ for every $n$ under which both the measure $\p_n$ and the function $f_n$ are invariant? Here is our answer for iid sequences in a probability space $(\Omega,\pi)$:

\bth[Clue of transitive functions]\label{t.cluegen}
Let $f\in L^2(\Omega^{V},\pi^{\otimes V})$, and suppose that there is a subgroup $\Gamma\leq \mathrm{Sym}_V$ acting on $V$ transitively such that $f$ is invariant under its action. If $U\subseteq V$, then
$$\clue(f\;|\;U)\leq \frac{\left|U\right|}{|V|}.$$  
In particular, sparse reconstruction for transitive functions of iid variables is not possible.
\eth

We will first give a proof for the case when $\pi$ is the uniform measure on $\{-1,1\}$, using the Fourier spectral sample, introduced in \cite{GPS}. We will then generalize this proof to general product measures using the Efron-Stein decomposition from \cite{ES} (see also \cite[Section 8.3]{OD}). Given these notions of a spectral sample, the proof turns out to be surprisingly simple (see Section~\ref{s.L2}). However, it is quite rigid, using transitivity in an essential way. Transitivity can be relaxed to quasi-transitivity if every orbit has a size comparable to the entire set (see Theorem~\ref{cluethm}), but, for instance, the condition that each variable has the same small influence is already not enough, as shown by the example of Remark~\ref{smallinflu} below. That example is even quasi-transitive, with two orbits, one being much larger than the other. We can fully avoid transitivity only if $U$ is not a carefully chosen deterministic set, but is random:

\bth[No reconstruction from sparse random sets]\label{t.cluerandom}
Let $f\in L^2(\Omega^{V},\pi^{\otimes V})$ be any function. Let $\mathcal{U}$ be a random subset of $V$, independent of the $\sigma$-algebra $\mathcal{F}_{V}$. Then 
$$
\E[\clue(f \; |\;\mathcal{U})] \leq  \delta(\mathcal{U}),
$$
where $\delta(\mathcal{U}) := \max_{j \in V}{\p[j \in \mathcal{U}]}$ is called the revealment of $\mathcal{U}$.
\eth

This notion of revealment was introduced in \cite{SS} for randomized algorithms computing Boolean functions by asking bits one-by-one, allowed to use the information already obtained in choosing which bit to ask next, along with extra randomness. Many interesting functions are known to have small revealment algorithms computing their values, and the key discovery of Schramm and Steif in \cite{SS} was that such functions are necessarily noise sensitive. Although this is not the usual definition, noise sensitivity can be defined in terms of $\clue$, and then the result of \cite{SS} can be stated as follows: if $f_n$ can be computed with a randomized algorithm with revealment $\delta_n$, and $\cB^{1-\eta_n}$ is an iid Bernoulli$(1-\eta_n)$ subset of $V_n$, with $\eta_n / \sqrt{\delta_n} \to \infty$, then $\E\big[ \clue(f_n \given \cB^{1-\eta_n}) \big] \to 0$. 

If $\cU$ was a small revealment subset, independent of $\cF_V$, which had a clue close to 1 about $f$, then asking the bits in $\cU$ would be a randomized algorithm that approximately computes $f$, hence the above theorem from \cite{SS} would say that $f$ is noise sensitive. But then, even the high density random set $\cB^{1-\eta}$ would have a small clue, so do we not get a contradiction to $\cU$ having large clue, obtaining a proof of Theorem~\ref{t.cluerandom} immediately from \cite{SS}? 

The answer is ``no'', for two reasons. One, getting from non-vanishing clue to a clue close to 1 does not seem to be an obvious matter (see Question~\ref{q.cpshape} at the end of the paper). Two, a small revealment random set $\cU$ might be quite different from an iid Bernoulli random set. A trivial example is when we ask all the bits with probability $\delta$ and none of the bits with probability $1-\delta$. This simple strategy achieves an average clue of $\delta$ for any function (the best possible according to Theorem \ref{t.cluerandom}), in contrast with asking an iid Bernoulli subset of the bits  with fixed density $\delta$, which has a vanishing average clue whenever the sequence of functions in question is noise sensitive.   

The small revealment theorem of \cite{SS} nevertheless suggests that if we want interesting non-transitive functions for which sparse reconstruction is possible, then we should probably look for noise-sensitive examples. A central example in the theory (see \cite{GS}) is left-to-right crossing in a box in critical planar percolation. Here a lot is known about the spectral sample \cite{GPS}, but still, the proof of Theorem~\ref{t.cluegen} does not generalize in a straightforward way. We will nevertheless answer Benjamini's question by proving in Section~\ref{s.perc} that there is no sparse reconstruction, with an argument that formalizes one's natural feeling that left-to-right crossing is not that far from being transitive, because the boundaries of the box should not play an important role. Namely, if one embeds an $n\times n$ box into a torus, and a left-to-right crossing occurs in the box, then the crossing is unlikely to stop exactly at the boundary, hence it is likely to work also in boxes shifted by a small macroscopic distance $\delta n$ on the torus. Thus, in a certain sense, the function is almost-quasi-transitive.  (See Proposition~\ref{l.almosttran} for a precise formulation of this notion.)
We then use Theorem~\ref{t.cluerandom} to conclude:

\bth[No sparse reconstruction in percolation]\label{t.clueperc}
There is no sparse reconstruction for $f_n$ being the indicator of left-to-right crossing in the box $[n]^2$ in critical bond percolation on $\Z^2$.
\eth

The same argument works for left-to-right crossing in the $n\times n$ rhombus $\mathsf{Rh}_n$ for critical site percolation on the triangular lattice. In that setting, it was proved in \cite{GPS} that if $U_n \subset \mathsf{Rh}_n$ is missing at least one site from every ball of radius $n^{3/8-\eps}$, for any fixed $\eps>0$, then its clue goes to 0. Consequently, if $|U_n| < n^{3/4-\eps}$, then its clue goes to 0. This is a very weak corollary, but it is unclear how to get a stronger result from the methods of \cite{GPS} in terms of just $|U_n|$. Furthermore, it was shown in \cite{PSSW}, using an inequality from \cite{ODS}, that even if we ask the bits in an adaptive way, then getting a clue close to 1 needs at least $n^{3/2+o(1)}$ queries. However, just as before, it is not clear to us if this applies also to getting any positive clue. 

If the revealment $\delta(\cU)$ is small, then the small expected clue provided by Theorem~\ref{t.cluerandom} implies that the clue is small with high probability. What happens if we only know that $\delta(\cU)$ is bounded away from 1? The famous ``It Ain’t Over Till It’s Over'' theorem of~\cite{MOO} says that, for sequences of functions with low maximal influence, for an arbitrary small (but fixed) $\eta$, the clue of a Bernoulli random subset $\cB^{1-\eta}$ is bounded away from 1 with high probability. However, for the $i$th dictator $\Dict_i(\omega):=\omega_i$, the expected clue $\E\big[\clue(\Dict_i \given \cB^{1-\eta})\big]=1-\eta$ is bounded away from 1, but the stronger result fails, since $\p\big[\clue(\Dict_i \given \cB^{1-\eta})=1\big]=1-\eta$.

\medskip
Natural as it may seem, $L^2$-clue is obviously not the only possible way to quantify the information content of a subset of coordinates about a function. For discrete random variables, entropy is a good alternative. 

\bde[I-clue] \label{d.Infoclue}
Let $X_V:=\{X_v :\;v \in V\}$ be a finite family of discrete random variables on some probability space $(\Omega^V,\,\p)$, and for some  $f: \Omega^V \lora \R$ consider the random variable $Z = f(X_V)$. The information theoretic clue (I-clue) of $f$ with respect to $U \subseteq V$  is
$$
\clue^{I}(f\;|\;U)= \frac{I(Z : X_U )}{I(Z : X_V )} = \frac{I(Z : X_U )}{H(Z)},
$$
where $H(Z):= - \sum_z \p[Z=z] \log \p[Z=z]$ is the entropy, and $I(Z : X) := H(Z)+H(X)-H(Z,X)$ is the mutual information.
\ede

The analogue of Theorem~\ref{t.cluegen} turns out to hold also in this setting. Of course, the proof now is not via spectral considerations --- we will use Shearer's entropy inequality (see \cite{Shearer} or \cite[Theorem 6.28]{LP}).

\bth[I-clue of transitive functions]\label{t.Infoclue}
Let $\left\{ X_v  : v \in V \right\}$ be discrete valued i.i.d.~random variables with finite entropy. Let $f: \Omega^V \longrightarrow \R$ be a transitive function and $Z = f(X_V)$. Then
\beq
\clue^{I}(f\;|\;U)\leq \frac{|U|}{|V|}.
\eeq
\eth

One can define further notions of clue, based on different distances between the distributions  $X_U$ and $X_U \;|\; Z$, such as total variation distance and Kullback-Leibler divergence, and analogous small clue theorems can be proved for iid variables. See Sections~\ref{s.otherclue} and~\ref{s.infoclue}. For typical non-degenerate random variables and functions $f$, the different notions of clue are comparable to each other, hence the different ``no sparse reconstruction'' theorems are more-or-less equivalent to each other. It may be nevertheless interesting to explore extreme examples, where sparse reconstruction is possible in one setting, but not in another (in particular, see Question~\ref{q.L2I}).
\medskip

It may seem very surprising that the exact same bound $|U|/|V|$ shows up for completely different notions of clue. The explanation is that there is a common generalization of these results, inspired by the Shapley value in cooperative game theory (see, e.g., \cite{BDT, PS}):

\bth[Clue from convex games] \label{t.coopclue} 
Let $v: 2^V \rightarrow [0,1] $ be a supermodular set function (cooperative game)  on the finite set $V$,  i.e., 
$$v(S)+v(T) \leq v(S \cup T)+v(S \cap T).$$
Let $\Gamma$ be a group acting on $V$ transitively, and suppose that the cooperative game $v$ is invariant under the action of $\Gamma$. 
Then, for any $U \subseteq V$,
$$
v(U) \leq  \frac{|U|}{{|V|}}v(V).
$$
\eth

If $\p$ is a measure on $\Omega^{V}$, we can define a cooperative game by $v_f(S) := \clue(f\;|\;S)$, with different notions of $\clue$. In particular, it is not difficult to verify that in case $\p$ is an iid measure, supermodularity is satisfied  when $v_f(S)$ is defined either as the $L^2$-clue or the $I$-clue.

The proof of Theorem~\ref{t.coopclue} is combinatorial, by induction, quite similar to the proof of Shearer's inequality. See Section~\ref{s.coop}.
\medskip

Our motivation for defining different notions of clue and proving the corresponding small clue theorems was not just abstract curiosity. In forthcoming work, we will study sparse reconstruction for non-iid measures $\p$ --- primarily the Ising model and factor of iid spin systems. As we have explained, a small clue theorem may be considered as a baby noise sensitivity result. However, discrete Fourier analysis breaks down for non-iid measures, hence it is highly desirable to come up with possible replacements. For instance, we will prove a small clue theorem for any high temperature Curie-Weiss model (the Ising model on the complete graph) using the I-clue and a version of Theorem~\ref{t.Infoclue}.
\medskip

\noindent{\bf Acknowledgments.} We are grateful to Itai Benjamini for his inspiring question that started this project. Our work was supported by the ERC Consolidator Grant 772466 ``NOISE''. During most of this work, PG was a PhD student at the Central European University, Budapest. Thanks to Christophe Garban and Bal\'azs Szegedy for reading the thesis, to Ohad Klein for a comment and a reference, and to two fantastic referees for their many corrections and suggestions.

\section{$L^2$-clue and sparse reconstruction for transitive functions}\label{s.L2}

\subsection{The Fourier-Walsh expansion and the Spectral Sample} \label{Fourierm}
We introduce a function transform on the hypercube which turns out to be an essential tool in the analysis of Boolean functions. We still consider the uniform measure $\p_{1/2}:=(\frac{1}{2}\delta_{-1} + \frac{1}{2}\delta_{1} )^{\otimes V}$. We can introduce the natural inner product $(f,g) = \E[fg]$ on the space of real functions on the hypercube.

\begin{definition} [Fourier-Walsh expansion]
For any $f \in L^2( \{-1,1\}^{V}, \p_{1/2}) \; \text{and} \; \omega \in\{-1,1\}^V$
\begin{equation}
f(\omega)=\sum_{S\subset V} \widehat{f}(S) \chi_S(\omega), \quad \quad \chi_S(\omega):=\prod_{i\in S}\omega_i\, \; (\text{and} \; \chi_S(\emptyset):= 1).
\end{equation}
\end {definition}

This is in fact the Fourier transform on an Abelian group: the event space is naturally identified with the group $\mathbb{Z}_{2}^{V} $ by assigning a generator $g_x$ to every $x \in V$, and the functions $\chi_S$ are the characters of $\mathbb{Z}_{2}^{V}$.

It is straightforward to check that the functions $\chi_S$ form an orthonormal basis with respect to the inner product, so Parseval's formula applies and therefore 
$$
 \sum_{S \subseteq V}{\widehat{f}(S)^2} =\left\|f\right\|^2 .
$$
Noting that $\widehat{f}(\emptyset) = \E[f]$, we also have 
\beq \label{VarFourier}
 \Var(f) = \sum_{\emptyset \neq S  \subseteq V}{\widehat{f}(S)^2}.  
\eeq 
For a subset $T \subseteq V$ let us denote by $\mathcal{F}_T$ the $\sigma$-algebra generated by the bits in $T$. So $\mathcal{F}_T$ expresses knowing the coordinates in $T$. It turns out that the conditional expectation  of any  function $f: \{-1,1\}^{n} \longrightarrow \R \;$  with respect to  $\mathcal{F}_T$ can be expressed in terms of the  Fourier-Walsh expansion; see \cite{GS}:
\beq \label{condexpFour}
\E[f \,|\, \mathcal{F}_T]= \sum_{S\subseteq T} {\widehat{f}(S)\chi_S}.
\eeq
The proof is fairly simple: we only need to observe that, if $S\subseteq T$, then $\E[\chi_S \,|\, \mathcal{F}_T] =\chi_S$, and in any other case  $\E[\chi_S \,|\, \mathcal{F}_T] =0$.

Using \eqref{VarFourier} we get a concise spectral expression for the variance of the conditional expectation:
\beq \label{VarFourierCond}
\Var(\E[f \,|\, \mathcal{F}_T])= \sum_{\emptyset \neq S\subseteq T} {\widehat{f}(S)^2}.
\eeq 

It turns out to be useful to think about the squared Fourier coefficients $\widehat{f}(S)^2$ as a measure on all the subsets of the coordinates. It is usually normalized  to get a probability measure. The random subset $\Spec_{f}$ distributed accordingly is called the spectral sample of $f$. 

\bde [Spectral sample]
Let $f \in L^2( \{-1,1\}^{V}, \p_{1/2})$. The spectral sample $\Spec_f$ of $f$ is a random subset of $V$ chosen according to the distribution 
$$
\mathbb{P}[\Spec_{f} = S] = \frac{\widehat{f}(S)^2}{\left\|f\right\|^2}, \quad \text{for any} \; S \subseteq V.
$$
\ede

The advantage of this concept is that it introduces a  rather compact language, where some important concepts  admit straightforward translations. 
The notion of clue, in particular, translates  well to the Spectral Sample language. Using \eqref{VarFourier} and  \eqref{VarFourierCond} we get that
\beq \label{clueFourier}
\clue(f\;|\;U) =\p [\Spec_{f} \subseteq U \,|\,\Spec_{f} \neq  \emptyset].
\eeq
This equation, as we shall see, is one of the key observations in the  proof of  Theorem \ref{cluethm}.

\subsection{No sparse reconstruction for transitive functions of fair coins}

The following theorem provides a sharp upper bound on the clue of not only Boolean, but general real-valued transitive functions. The proof is surprisingly short and it demonstrates the power of the notion of spectral sample.

\bth [Clue of transitive functions] \label {cluethm}
Let $V$ be a finite set and $\Gamma \curvearrowright V$ a group acting on $V$  and for a $v \in V$ let $\Gamma \cdot v = \{ v^{\gamma} : \gamma \in \Gamma \}$ be the orbit of $v$. Let $f: \{-1,1\}^{V} \longrightarrow \R$ be invariant under the action of $\Gamma$. Then for any $U\subseteq V$
$$ 
\clue(f \given U)\leq \ \frac{\left|U\right|}{\min_{v \in V}{|\Gamma \cdot v|}}.
$$ 
In particular, if  $f$ is transitive,  then
\beq \label{cluettrans}
\clue(f\;|\;U)\leq \frac{\left|U\right|}{\left|V\right|}.
\eeq
\eth

\begin{proof}
Let $X$ be a uniformly random element from the spectral sample $\Spec_f$ of $f$ conditioned on being non-empty. Because $f$ is invariant under the action of $\Gamma$ by assumption, if $v$ and $w$ are both in the same orbit $\Gamma \cdot u$  then using that the distribution of $X$ is also invariant with respect to $\Gamma$ we have  
$$
\Tilde{\p}[X=v] = \Tilde{\p}[X=w],
$$
where $\Tilde{\p}$ denotes the probability measure conditioned on $\{ \Spec \neq  \emptyset\}$.
Therefore
$$
\sum_{s \in  \Gamma \cdot u} {\Tilde{\p}[X=s]} = | \Gamma \cdot u|\Tilde{\p}[X=u] \leq 1,
$$
which gives the bound $\Tilde{\p}[X=u] \leq 1/ | \Gamma \cdot u|$.

Using \eqref{clueFourier}  we get the following:
\beq 
 \clue(f\;|\;U) =\p [\Spec \subseteq U \,|\,\Spec \neq  \emptyset] \leq \Tilde{\p}[X \in U]  = \sum_{u\in U}\Tilde{\p}[X=u] \leq \frac{\left|U\right|}{\min_{v \in V}{\left|\Gamma \cdot v\right|}}.
\eeq
Now \eqref{cluettrans} is an obvious consequence using that when $\Gamma$ is transitive then  $\Gamma \cdot v = V$ for all $v \in V$. 
\end{proof}

\begin{rem}
The bound in Theorem \ref{cluethm} is sharp, as it is testified by the function $\sum_{v \in V}{\omega_v}$.
\end{rem}

It is worth to emphasize that the result does not only apply for sequences of Boolean functions, but also for any sequences of real-valued functions, no matter bounded or not.

\begin{rem}\label{smallinflu}
There is no obvious way to relax the condition of transitivity. Let $f : \{-1,1\}^{V} \longrightarrow  \{-1,1\} $ and $j \in V$. The influence of the coordinate $j$ is $I_j(f) :=\p[f(\omega) \neq f(\omega^{j}) ]$, where $\omega^{j}$ is the same as $\omega$ except its $j$th coordinate is flipped. 

We now sketch an example of a sequence of Boolean functions where the  influences $I_{j}({f_n})$ are (almost) equal for every $n$, however there is a sparse subset of coordinates $U_n$ (i.e., $\lim_n {\frac{|U_n|}{|V_n| }}=0$) such that $ \lim_{n}\clue (f_n \given U_n) = 1$.

Let $a_n$ be a sequence of integers such that $a_n \rightarrow \infty$. Let us define the asymmetric majority functions
$$
\Maj^{a_n}_n(\omega) =
\left\{
	\begin{array}{ll}
1	 & \mbox{if } \sum_{i=1}^n{\omega(i)} > a_n\sqrt{n} \\
		-1 & \mbox{if } \sum_{i=1}^n{\omega(i)} < a_n\sqrt{n}.
	\end{array}
\right.
$$
One can choose $a_n$ in such a way that 
$$
I_{i}(\Maj^{a_n}_{n})=\frac{{{n} \choose {n/2+2a_n\sqrt{n}}}}{2^n} \sim \frac{1}{n^{2/3}}  
$$
holds. 

Furthermore, define the Boolean function $\Tr^{l,k}: \left\{ -1,1 \right\}^{lk} \longrightarrow  \left\{0, 1 \right\}$ as follows: we group the bits in $k$ $l$-element subsets, the so called tribes. The function takes on $1$ if there is a tribe $T$ such that $\omega(i)=1$ for every $i \in T$, and $-1$ otherwise.

$\Tr^{l_n, k_n}$  is known to be balanced if $l_n = \log n - \log \log n$ and $k_n = n/l_n$. Let us denote this balanced version of the tribes on $n$ bits by $\Tr_n$. An straightforward calculation shows that $I_{i}(\Tr_n) \sim \frac{\log n}{n}$.  

Take a disjoint union $V_n = M_n \cup T_n$, with $|M_n| = m_n$ and  $|T_n| = t_n$. Now we define our function as follows: 
$$
f_n :=
\left\{
	\begin{array}{ll}
\Maj^{a_n}(\omega_{M_n})	 & \mbox{if } \Tr(\omega_{T_n}) = 1 \\
		\Maj^{-a_n}(\omega_{M_n}) & \mbox{if } \Tr(\omega_{T_n}) = -1.
	\end{array}
\right.
$$
We adjust the size of $M_n$ and $T_n$ in such a way that the influence of each coordinate is the same. So we have the equation $\frac{\log t_n}{t_n} = \frac{1}{m_n^{2/3}}$, or equivalently
$$
m_n = \left(\frac{t_n}{\log t_n}\right)^{3/2}.
$$
So the density of $T_n$ goes to $0$ compared to $|V_n| = t_n + m_n$. At the same time, from the Central Limit Theorem it is clear that $\lim_{n}\p[\Maj^{a_n}_{m_n} = 1] = 0$ and $\lim_{n}\p[\Maj^{-a_n}_{m_n} = 1] = 1$. Consequently, $\lim_{n} \clue ( f_n \given T_n ) = 1$.  
\end{rem}

\begin{rem}
We point out an interpretation of the random element $X$ of the spectral sample appearing in the proof of Theorem \ref{cluethm}. This setup also has some interesting connections with one of the key lemmas in Chatterjee's book on superconcentration and chaos \cite{Ch}. 

For a function $f  : \{-1,1\}^{V} \longrightarrow \R$ we define the stability of $f$ at level $p$ as 
$$
\Stab_f(p) := \sum_{\emptyset \neq S\subseteq V}{\widehat{f}(S)^2 p^{\left|S\right|}}. 
$$
(This is a small modification of the definition in \cite{OD}.) Let us denote by $\omega^{1-p}$ the random vector which we obtain from $\omega$ by resampling each of its bits independently with probability $1-p$. With this notation,  we clearly have $\Stab_f(p) = \Cov (f(\omega),f(\omega^{1-p}))$.

At the same time, it is also the expected clue of a Bernoulli random set of coordinates $\mathcal{B}^{p}$ of density $p$:  $\frac{\Stab_f(p)}{\Var(f)} = \E [ \clue (f \;|\; \mathcal{B}^{p})]$.

Stability can be generalized as a polynomial of $|V|$ variables. Then the quantity 
$$
\frac{\Stab_f(\mathbf{x}) }{\Var(f)} = \frac{1}{\Var(f)}\sum_{\emptyset \neq S\subseteq V}{\widehat{f}(S)^2 \prod_{i\in S}x_i}
$$ 
can be interpreted as the expected clue of a random subset where the bit $i$ is selected with probability $x_i$, independently from other bits.

Denote by $\overline{p}$ the vector with all coordinates equal to $p$ and for a $j \in V$ take the partial derivative of $\Stab_f(\overline{p})$ with respect to the $j$th coordinate. We obtain that 
$$
\frac{\partial\Stab_f(\overline{p})}{\partial p_j } =  \sum_{S \ni j} {\widehat{f}(S)^2p^{|S| -1} }.
$$
Now here is the relationship with $X$, the uniformly random element of the spectral sample:
\beq \label{XStab}
\int_{0}^{1} {\frac{\partial\Stab_f(\overline{p})}{\partial p_j } dp}= \sum_{ S \ni j}{\widehat{f}(S)^2\frac{1}{|S|}} = \Var(f) \p[X=j]. 
\eeq
The above quantity can be understood as the average increase in clue over all $p$ values, induced by  a small increase in the probability of selecting $j$ into the random set. This interpretation becomes even more explicit in the cooperative game theory framework (see Proposition \ref{Xisshapley} below).

Now we get to the connection with Chatterjee's work. Let $f,g  : \{-1,1\}^{V} \longrightarrow \{-1,1\}$ be monotone Boolean functions. We say that a coordinate $j$ is pivotal for $f$ (given $\omega$) if switching the coordinate $j$ changes the value of $f$. The set  $\Piv_f(\omega)$ of  all pivotal coordinates  is called the pivotal set.  We start by expressing $\p[ j \in \Piv_f(\omega) \cap \Piv_g(\omega^{1-p})]$ in terms of the Fourier-Walsh transform.

Observe that for any monotone $f: \{-1,1\}^{V} \longrightarrow \{-1,1\}$, we have 
$$
\nabla_j f(\omega) := f(\omega| \omega_j =1 ) - f(\omega| \omega_j = -1 ) = \sum_{S \ni j} {\widehat{f}(S) \chi_{S \setminus j}(\omega)}.
$$
As $j$ is in $\Piv_f(\omega)$ if and only if $\nabla_j f(\omega) =2$, and otherwise $\nabla_j f(\omega) =0$, we get that 
$$
\1_{j \in \Piv_f(\omega)}  = \frac{1}{2} \sum_{S \ni j} {\widehat{f}(S) \chi_{S \setminus j}(\omega)}.
$$
Now recall that 
$$
\E[ \chi_{T}(\omega)\chi_{S}(\omega^{1-p})]  =
\left\{
	\begin{array}{ll}
0	 & \mbox{if } T \neq S, \\
		p^{|S|} & \mbox{if }  T = S,
	\end{array}
\right.
$$   
and thus, whenever $f$ and $g$ are monotone, we have
\beq \label{pivcorrfour}
\p[ j \in \Piv_f(\omega) \cap \Piv_g(\omega^{1-p})] = \E[\1_{j \in \Piv_f(\omega)} \1_{ j \in  \Piv_g(\omega^{1-p})} ] =  \frac{1}{4} \sum_{S \ni j} {\widehat{f}(S)\widehat{g}(S) p^{|S| -1} }.
\eeq
(We note that this formula is almost a generalization of Lemma 2.7 in \cite{RS}.)
Using that $\sum_{j \in V} {\p[ j \in \Piv_f(\omega) \cap \Piv_g(\omega^{1-p})]} =\E[  |\Piv_f(\omega) \cap \Piv_g(\omega^{1-p})|]$, we get from \eqref{pivcorrfour} that
\beq \label{covlemma}
\int_{0}^{1} {\E[  |\Piv_f(\omega) \cap \Piv_g(\omega^{1-p})|] dp}= \frac{1}{4} \sum_{j \in V}  \left( \sum_{ S \ni j}\widehat{f}(S)\widehat{g}(S) \frac{1}{|S|}\right) =\frac{1}{4}\Cov(f,g). 
\eeq
This is essentially a special case of Lemma 2.1 from \cite{Ch} (referred to as ``covariance lemma''), where the Markov process is the random walk on the hypercube. At the same time, setting $g = f$, by \eqref{XStab} we have
$$
\p[X = j] = \frac{1}{\Var(f)}\int_{0}^{1} {\frac{\partial\Stab_f(\overline{p})}{\partial p_j } dp}=\frac{4}{\Var(f)} \int_{0}^{1} {\p[ j \in \Piv_f(\omega) \cap \Piv_f(\omega^{1-p})] dp},
$$
which is a coordinate-wise localized version of the covariance lemma.

We also note that the threshold saddle vertex introduced in \cite{Riv} for level set percolation is a random pivotal vertex, so it might be considered as a real (non-Fourier) space analogue of our $X$, and it is also closely related to Chatterjee's covariance lemma.
\end{rem}

\subsection{No sparse reconstruction in general product measures}\label{ss.genprod}

One may ask whether a result similar to Theorem \ref{clue} can be derived in case we replace $ \{-1,1\}$  in the domain with another space, or if we replace the product measure with some other measure. A natural idea in this direction is to try to generalize the concept of spectral sample. We might take again equation~\eqref{clueFourier} as a starting point.

In the previous section we denoted the binary coordinates by $\omega_i$ for some $i \in V$. In the general setting, when the coordinates are general real-valed random variables, we shall use $X_i$ instead.       

Observe that the quantity $\clue(f|U)$  is well defined for any $U \subseteq V$ on any product space $X^{V}$, no matter what the probability measure is. So one could try to use equation \eqref{clueFourier} as the definition for a generalised spectral sample.  As the probabilities $\p[\Spec \subseteq U]$ are known for all $U$, one can also calculate the probabilities $\p[\Spec = T]$ for all $T$. Once we have this generalised spectral sample in hand (depending on the function, the space and the underlying measure) we might be able to repeat the argument in the proof of Theorem \ref{cluethm}.

Unfortunately this strategy fails in general. The problem is that nothing guarantees that the quantities $\p[\Spec = T]$ that we get from the M\"obius inversion are non-negative. Nevertheless, in case the underlying measure is a product measure, the above strategy works as the quantities  $\p[\Spec = T]$ turn out to be non-negative. As we will show, this follows directly from the so-called Efron-Stein decomposition (see \cite{ES} or \cite[Section 8.2]{OD}), a generalization of the Fourier-Walsh transform for general product measures. 

We will need the following simple observation, true only for product measures, which turns out to be crucial. In fact, as we shall see, the Efron-Stein decomposition, as well as the possibility of a spectral sample, ultimately depend on Fubini's Theorem.

\bl\label{keystep} 
Let $f\in L^2(\Omega^n,\pi^{\otimes n})$ and let $K,L \subseteq [n]$. Then
$$
\E[\E[f \,|\, \mathcal{F}_L] \,|\, \mathcal{F}_K] = \E[f \,|\, \mathcal{F}_{L \cap K}].
$$
\el
\bpf
Rewriting the conditional expectations as integrals, and using Fubini's theorem,
\[
\int_{X^{K^{c}}}{ \left( \int_{X^{L^{c}}}{f(X_L, x_{L^{c}} dx_{L^{c}}) } \right) dx_{K^{c}} } = \int_{X^{K^{c} \cup L^{c}}}{f(X_{L \cap K}, x_{K^{c} \cup L^{c}}) dx_{K^{c} \cup L^{c}}  }. 
\qedhere
\]
\epf

\bth [Efron-Stein decomposition \cite{ES}]\label{EfronStein}
For any $f\in L^2(\Omega^n,\pi^{\otimes n})$, there is a unique decomposition
$$f=\sum_{S\subseteq [n]} f^{S}\,,$$
where $f^{S}$ is a function that depends only on the coordinates in $S$, and if $S \not\subseteq T$ then $(f^{S},g)=0$ for any $\mathcal{F}_T$--measurable function $g$.
\eth

For completeness, we include a proof, following the ideas from \cite{OD}, but with our notation, and pointing out the key role of Lemma~\ref{keystep}.

\bpf
Observe that, assuming that the sought decomposition exists, we have  (as in the case of the hypercube, see \eqref{condexpFour})
$$
\E[f \,|\, \mathcal{F}_T] = \sum_{S \subseteq T}{f^{S}}.
$$
Indeed, this is clear from the fact that $\E[f^{S} \,|\, \mathcal{F}_T]$ is $f^{S}$ in case $S \subseteq T$ (as in this case $f^{S}$ is $\mathcal{F}_T$--measurable) and $0$ otherwise (because of the  second property of the decomposition). 

We can use this to find the functions $f^{S}$ via M\"obius inversion (in this case, an exclusion-inclusion principle) from the conditional expectations. So, define $f^S$ in the only possible way (already showing uniqueness), by
$$
f^{S} := \sum_{L \subseteq S}{ (-1)^{|S|-|L|}\E[f \,|\, \mathcal{F}_L]  }.
$$
It is obvious from the construction that $f^{S}$ only depends on coordinates in $S$. We now show that if $g$ is $\mathcal{F}_T$--measurable and $S\setminus T \neq \emptyset$ then $f^{S}$ and $g$ are orthogonal. 

Pick an $i \in S\setminus T$ and write the above inner product as
$$
\E[g f^{S}] = \sum_{L \subseteq S\setminus \left\{ i \right\} }{ (-1)^{|S|-|L|} \left( \E[g\E[f \,|\, \mathcal{F}_L]]- \E[g\E[f \,|\, \mathcal{F}_{L \cup \left\{ i \right\}}]] \right)},
$$
using that $(-1)^{|S|-|L|}$ and $(-1)^{|S|-|L\cup \left\{ i \right\}|}$ have opposite signs.

Now we show that $\E[g\E[f \,|\, \mathcal{F}_L]] = \E[g\E[f \,|\, \mathcal{F}_{L \cup \left\{ i \right\}}]$ which implies that $\E[g f^{S}]=0$ and thus concludes the proof.

First note that
$$
\E[g\E[f \,|\, \mathcal{F}_L]] = \E[ \E[g \E[f \,|\, \mathcal{F}_L]\,|\, \mathcal{F}_T] ] =  \E[g \E[f \,|\, \mathcal{F}_{T \cap L}]], 
$$
using that $g$ is $\mathcal{F}_T$--measurable and by Lemma \ref{keystep}. 

As $T \cap (L \cup \left\{ i \right\}) = T \cap L$ (recall that  $i \notin L$ and $i \notin T$), we can conclude:
\begin{align*}
\E[g \E[f \,|\, \mathcal{F}_{T \cap L}]]& = \E[g \E[f \,|\, \mathcal{F}_{T \cap (L \cup \left\{ i \right\})}]] \\
&= \E[ \E[g \E[f \,|\, \mathcal{F}_{ L \cup \left\{ i \right\} }]\,|\, \mathcal{F}_T] ] = \E[g\E[f \,|\, \mathcal{F}_{L \cup \left\{ i \right\}}]]. \qedhere
\end{align*}

\epf
Observe that this is indeed a generalization of the Fourier-Walsh transform, with $f^{S} = \widehat{f}(S) \chi_S$. What is important for our purpose is that we can again define a Spectral Sample  $\p[\Spec = S]:= \frac{\|f^{S}\|^2}{\|f\|^2}$ for every square-integrable function, as in the case of the hypercube and thus Theorem \ref{cluethm} generalizes for product measures. 

\bth [Small clue theorem for general product measures] \label{cluethm2}
Let $f\in L^2(\Omega^n,\pi^{\otimes n})$ and suppose that there is a $\Gamma \leq S_n$ acting on the $n$ copies of $\Omega$  such that $f$ is invariant under this action. If $U\subseteq [n]$, then
$$ 
\clue(f\;|\;U)\leq \ \frac{\left|U\right|}{\min_{1 \leq j \leq n }{|\Gamma \cdot j|}}.
$$ 
In particular, if  $f$ is transitive,  then
$$  \clue(f\;|\;U)\leq \frac{\left|U\right|}{n}.  $$  
\eth
The proof is exactly the same as for Theorem \ref{cluethm}, the only difference being that we need to use the Efron-Stein decomposition instead of the Fourier-Walsh transform to build the spectral sample. 
\medskip

We close this section by giving a generalization of Theorem~\ref{cluethm2} that will play a key role in the proof of Theorem \ref{t.clueperc}. One of the advantages of this version is that it avoids the notion of transitivity altogether. In this setup we consider the average clue with respect to a random subset of coordinates.  It is important that the subset is sampled independently from the coordinate values. The size of a random subset $\mathcal{U}$ is measured in revealment, that is $\delta(\mathcal{U}) := \max_{j \in V}{\p[j \in \mathcal{U}]}$, a concept introduced in \cite{SS} for randomized algorithms. Our goal is to find for a given function a random subset of coordinates with fixed revealment which has as high expected clue as possible.

For any function $f$ one may use the random subset which is the entire coordinate set $V$ with probability $\delta$ and $\emptyset$ with probability $1-\delta$. The revealment for this trivial strategy is  $\delta$ as well as the expected clue we achieve. Somewhat surprisingly, it turns out that this is the best we can do.

\bth[No reconstruction from sparse random sets] \label{NoRSR}
Let $f\in L^2(\Omega^{V},\pi^{\otimes V})$ be any function. Let $\mathcal{U}$ be a random subset of $V$, independent of the $\sigma$-algebra associated with $\pi^{\otimes V}$. Then 
$$
\E[\clue(f \; |\;\mathcal{F}_{\mathcal{U}})] \leq  \delta(\mathcal{U}),
$$
where $\delta(\mathcal{U}) := \max_{j \in V}{\p[j \in \mathcal{U}]}$ is called the revealment of $\mathcal{U}$ and $\mathcal{F}_{\mathcal{U}}$ denotes the smallest $\sigma$-algebra such that the set $\mathcal{U}$ and the random variables $X_j$, for all $j \in \mathcal{U}$ are measurable w.r.t.  $\mathcal{F}_{\mathcal{U}}$.
\eth

\bpf
The proof basically repeats the proofs of Theorems~\ref{cluethm} and~\ref{cluethm2}. Generate the Efron-Stein spectral sample $\Spec_f$ independently of $\mathcal{U}$, and let $X$ be a uniformly random element from $\Spec_f$ conditioned on being non-empty; $\Tilde{\p}$ denotes the respective conditional probability measure.  From the Efron-Stein analogue of \eqref{clueFourier}, it is easy to verify that $\E[\clue(f \; |\;\mathcal{U})] = \p [\Spec \subseteq \mathcal{U} \,|\,\Spec \neq  \emptyset]$. Therefore, using that $\mathcal{U}$ is independent of the $\sigma$-algebra of $\pi^{\otimes V}$, we get
\begin{align*}
  \E[\clue(f \; |\;\mathcal{U})] &\leq  \Tilde{\p}[X \in\mathcal{U}]   =\sum_{j\in [n]}\Tilde{\p}[X=j,\; j \in \mathcal{U}] \\
  &=\sum_{j\in [n]}{\Tilde{\p}[X=j] \p[j \in \mathcal{U}]} \leq\delta(\mathcal{U}) \sum_{j\in [n]} {\Tilde{\p}[X=j]} = \delta(\cU). \qedhere 
\end{align*}
\epf

This statement can be read in such a way that for product measures there is no reconstruction for any sequence of functions from a sparse sequence of random subsets of coordinates $\mathcal{U}_n$ (that is, for which $\delta(\mathcal{U}_n) \rightarrow 0$).

Furthermore, it is not difficult to see that Theorem~\ref{NoRSR} implies Theorem~\ref{cluethm2}. If $f$ is a transitive function, and $U \subseteq V$ is a fixed subset, then $\clue(f \; |\;U) =\clue(f \; |\;U^{\gamma})$  for any translation $\gamma \in \Gamma$. One can easily verify that a uniform random translate of $U$ is a random  subset with revealment $|U/|V|$, and thus we recover the bound in Theorem~\ref{cluethm2}. For the general, quasi-transitive case, we again consider a uniform random $\Gamma$-translate $\cU$ of $U$. If the orbits of $\Gamma$ on $V$ are denoted by $O_1,\dots,O_k$, and we let $U_i:=O_i\cap U$, then for any $v\in O_i$ we have $\p[v\in\cU]=|U_i|/|O_i|$. That is, the revealment of $\cU$ is at most $|U| / \min_i |O_i|$, and thus Theorem~\ref{NoRSR} gives us the general bound of Theorem~\ref{cluethm2}.

\subsection{Clue of almost transitive functions} 

In this section we present a few results that are necessary for the proof of our percolation result Theorem~\ref{t.clueperc}, but also have some interest of their own. The main topic is to extend Theorem~\ref{t.cluegen} by relaxing the condition that the function $f$ is invariant under the group action $\Gamma$. 

The following  lemma  basically states that in case two functions are highly correlated and one of them has high clue with respect to a $\sigma$-algebra, then the other function also has high clue  with respect to the same $\sigma$-algebra.

It is worth pointing out that the proof uses only basic facts from linear algebra. 
The geometric intuition is that in case the angle between two vectors is small, and a projection (i.e., the conditional expectation) does not decrease the norm of the first one too much, then it cannot  decrease much the norm of the other vector, either. We emphasize that this result holds in  general, that is, the underlying measure does not need to be a product measure.

\bl \label{proj}
Let $f, g \in  L^2(\Omega^V,\p)$ with 
$$
\Corr(f,g)\geq 1 -\eps. $$
Let $U \subseteq V$. If  
$$
\clue(f\;|\;U) \geq c, $$
then
$$
\clue(g\;|\;U)\geq  c  - 5\sqrt{\eps}.
$$ 

\el

\bpf
Without loss of generality we may assume that $\E[f]=\E[g]=0 $  and $\Var (f) = \Var (g) = 1$, and therefore we may use $\|\; \; \|^2$ instead of variance.

Using that the conditional expectation is an orthogonal projection, we have
\beq \label{pythag}
\|\E[f \,|\, \mathcal{F}_U] \|^2 + \|f -\E[f \,|\, \mathcal{F}_U] \|^2 = \| f\|^2,
\eeq
and, for every $\cF_U$-measurable $h$,
$$
\|f -\E[f \,|\, \mathcal{F}_U] \|^2 \leq \| f-h\|^2.
$$
Therefore, with the triangle inequality we get
\beq \label{triangle}
\|g -\E[g \,|\, \mathcal{F}_U] \| \leq \| g -\E[f \,|\, \mathcal{F}_U]\| \leq \| g - f \| + \| f -\E[f \,|\, \mathcal{F}_U]\|. 
\eeq
Now,
\begin{align*}
\|f - g \|^2 = \Var (f -g)  &=\Var (f) + \Var (g) - 2 \sqrt{\Var (f)\Var (g)}\Corr (f, g)  \\
&=2 (1-  \Corr (f, g))  \leq 2 \eps. 
\end{align*}
By assumption,
$$
\clue(f\;|\;U) = \frac{\|\E[f \,|\, \mathcal{F}_U] \|^2 }{ \| f\|^2} =\|\E[f \,|\, \mathcal{F}_U] \|^2 \geq c,
$$
so  \eqref{pythag} translates to the following bound: 
$$
\| f -\E[f \,|\, \mathcal{F}_U]\|^2 = \| f\|^2 - \|\E[f \,|\, \mathcal{F}_U] \|^2   \leq 1 - c.
$$
Plugging the estimates into \eqref{triangle} we can write (using that dividing by $\| g\|^2=1$ does not change the equation)
$$
\frac{\|g \|^2 - \|\E[g \,|\, \mathcal{F}_U] \|^2}{\| g\|^2} = \frac{\|g -\E[g \,|\, \mathcal{F}_U] \|^2}{\| g\|^2} \leq \left(\sqrt{2 \eps} +\sqrt{1-c} \right)^2, 
$$
and thus we get
$$
1- \clue(g\;|\;U) \leq 2\eps + 1-c + 2 \sqrt{2 \eps(1-c)} \leq  1-c + 2 \eps + 2 \sqrt{2 \eps},
$$
from which  it   is immediate (assuming $\eps \leq 1$) that 
\[
\clue(g\;|\;U)\geq  c -( 2 + 2 \sqrt{2})\sqrt{\eps} \geq c - 5\sqrt{\eps} . 
\qedhere
\]
\epf

Using Lemma~\ref{proj}, we can relax the condition of quasi-transitivity in Theorem~\ref{cluethm} to a certain almost-quasi-transitivity of the functions $f_n$. Namely, it is enough that there is a quasi-transitive group action $\Gamma_n \curvearrowright V_n$ with large orbits only, such that the translated functions 
$f_n^\gamma(\omega):=f_n(\omega^{\gamma^{-1}})$, for all $\gamma\in\Gamma_n$, fall into a bounded number of equivalence classes: if $f_n$ and $f_n^\gamma$ are equivalent, then their correlation is close to 1 (as a relaxation of being the same function, as in true transitivity).

\bpr\label{l.almosttran}
Let $V$ be a finite set and $\Gamma \curvearrowright V$ a group acting on $V$. 
Let $D \subseteq \Gamma$ and $L \subseteq \Gamma$ such that
\begin{enumerate}
\item $D = D^{-1}$ and $L = L^{-1}$,
\item $D \cdot L = \{ d \cdot \ell \;:\; d \in  D,  \; \ell \in L\} = \Gamma$.
\end{enumerate}
Suppose that $f: \{-1,1\}^{V} \longrightarrow \R$ has the  property that, for any $d \in  D$,
\beq \label{highcorrcond}
\Corr(f, f ^{d}) > 1- \eps.
\eeq
Then, for any $U \subseteq  V$,
\beq \label{e.LUclue}
\clue(f\;|\;U)\leq \ \frac{\left|L\right|\left|U\right|}{\min_{v \in V}{|\Gamma \cdot v|}} +5\sqrt{\eps}.
\eeq
In particular, if $\Gamma_n \curvearrowright V_n$ is a sequence of group actions such that $|V_n| / \min_{v\in V_n}|\Gamma_n\cdot v|$ is bounded, and $f_n:\{-1,1\}^{V_n} \longrightarrow \R$ is such that for every $\eps>0$ there exist subsets $D_{\eps,n}$ and $L_{\eps,n}$ that satisfy the above conditions 1., 2., and~(\ref{highcorrcond}), and $|L_{\eps,n}|$ remains bounded as $n\to\infty$, then there is no sparse reconstruction for $f_n$.
\epr

\bpf
Let $c:=\clue(f\;|\;U)$. An application of Lemma \ref{proj} for $f$ and 
$f ^{d}$ gives that, whenever  $d \in  D$,
$$
\clue(f ^{d}\;|\; U)\geq  c  - 5\sqrt{\eps}.
$$ 
As $\clue(f^{d}\;|\;U)=\clue(f\;|\;U^{d^{-1}}) = c$, we conclude that
$$
\clue(f\;|\;U^{d})\geq  c  - 5\sqrt{\eps},
$$
for every $d \in  D$.

Now let $W := \cup_{\ell \in L} {U^{\ell}}$. Take an arbitrary $\gamma \in \Gamma$. By our assumptions, it can be written as $\ell^\ast\cdot d^\ast $ for some $\ell^\ast \in  L$ and $d^\ast \in  D$. So clearly  we have 
$$
(U^{({\ell^\ast}^{-1})})^{\gamma} = U^{({\ell^\ast}^{-1})\ell^\ast d^\ast} = U^{d^\ast}.
$$
Using that ${\ell^\ast}^{-1}$ is also in $L$, this shows that
$U^{d^\ast} \subseteq  W^{\gamma}$.  Therefore,
$$
\clue(f\;|\;W^{\gamma}) \geq \clue(f\;|\; U^{d^\ast})\geq  c  - 5\sqrt{\eps}.
$$

Let $\mathcal{W}$ be a uniform $\Gamma$-translate of $W$. As in the last paragraph of Subsection~\ref{ss.genprod}, the revealment of $\mathcal{W}$ is at most $|W| / \min_{v\in V} |\Gamma\cdot v|$, and thus Theorem~\ref{NoRSR} tells us that
$$
c  - 5\sqrt{\eps} \leq \E[\clue(f\;|\;\mathcal{W})] \leq \frac{|L| |U|}{\min_{v\in V} |\Gamma\cdot v|},
$$
proving (\ref{e.LUclue}).

The statement of no sparse reconstruction follows from~(\ref{e.LUclue}), since if $|U_n|/|V_n|\to 0$, then both terms in the upper bound can be made arbitrarily small by first choosing $\eps$ small, then $n$ large enough.
\epf

\begin{rem}
The extension from transitivity to almost-quasi-transitivity does not use that the measure is a product measure: if the measure satisfies Theorem~\ref{NoRSR} in the sense that a vanishing revealment implies vanishing clue for transitive functions, then this remains true for almost-quasi-transitive functions.
\end{rem}

\section{Other approaches to measuring ``clue''}\label{s.otherclue}

\subsection{Significance and influence of subsets}\label{siginf}

We would like to make a small detour to discuss some possible alternatives to ``clue'' as defined in Definition \ref{clue}. Given  a Boolean function $f: \{-1,1\}^{V} \longrightarrow \{-1,1\}$ and an underlying probability measure $\p$, we want to quantify the amount of information a subset of the coordinates gives us about the function $f$. We will denote the size of the coordinate set $V$ by $n$.

We start with a sort of dual to $\clue$. 

\bde
The significance of a subset $U \subseteq V$ is
$$
\sig(f \;|\; U) = \frac{\E[\Var(f  \;|\;  \mathcal{F}_{U^c})]}{\Var(f)}
$$
\ede

We call it a dual because we have $\sig(f \;|\; U) = 1 - \clue(f \;|\; U^c)$. It expresses how much information we are still missing on average if we know the values of the bits outside of $U$. We have the following description of $\sig(f \;|\; U)$ in terms of the spectral sample:
$$
\sig(f \;|\; U) = \p[\Spec_f \cap U \neq \emptyset \;|\;  \Spec_f \neq \emptyset ].
$$
This shows that for product measures $\sig(f \;|\; U) \geq \clue(f \;|\; U)$. In general, this inequality does not hold: $\sig(f \;|\; U) < \clue(f \;|\; U)$ whenever $\clue(f \;|\; U) + \clue(f \;|\; U^c)>1$, which can easily happen if the underlying measure has lots of dependencies.  Also, Theorem \ref{cluethm} is not true if we replace $\clue$ by $\sig$. For example,  any subset $U \subseteq V$ has significance 1 with respect to the parity function $\chi_{V}$, which is obviously transitive.

We mention a similar concept introduced in \cite{BL}. For a subset $U \subseteq V$ the influence of $U$ is defined as follows: 
$$
I(f \;|\;U) = \p[f \;\text{is not determined by the bits on} \; U^{c}]
$$
Influence is much weaker than $\sig$, in the sense that it is  easier to have high influence than to have high significance: it is clear from the definition that, for any underlying measure, $I(f \;|\;U) \geq \sig(f \;|\;U)$. Like in social choice theory, one may think about coordinates as individual agents trying to influence the value (outcome) of $f$ by the values of the respective bits. In this framework the influence of a subset quantifies the probability that the set of agents in $U$ can change the value of $f$ by coordinating their values. While in this setting coordinates are allowed to cooperate, significance rather quantifies the average gain of information (measured in variance) for a uniformly random configuration of $U$.

We can again take the dual concept of influence, the combinatorial equivalent of clue, which is the probability that the subset $U$ is a witness. For a Boolean function $f: \{-1,1\}^{V} \longrightarrow \{-1,1\}$ and a configuration $\omega \in \{-1,1\}^{V}$ a subset $W \subseteq V$ is a witness for $f$ if $\omega_{W}$ already decides the value of $f$.
$$
W(f \;|\;U) = 1 - I(f \;|\;U^{c}) = \p[f \;\text{is determined by the bits on} \; U]
$$
Since the influence dominates significance, we have $W(f \;|\;U) \leq \clue(f \;|\; U)$  . But even for product measures, $I(f \;|\;U) \geq W(f \;|\;U)$ fails to hold (unlike the $\sig \geq \clue$ inequality, see above). For example, if $f(\omega) = \omega_1 \lor \omega_2 \lor \omega_3 $ and $U = \{ 1,2\}$, we have $I(f \;|\;U) < W(f \;|\;U)$. On the other hand, when $f = \chi_{\{ 1,2\}}$ and $  U = \{1\}$, clearly $I(f \;|\;U) > W(f \;|\;U)$.

There are still many questions to be investigated. For the left-right crossing event $\Cr_n$  for critical planar percolation, when $U_n$ is a sub-square, it is proved in \cite{GPS} that $I(\Cr_n \;|\;U_n) \asymp \sig(\Cr_n \;|\;U_n)$. For $\Maj_n$ on the other hand, this is not the case. As is easy to check, $I(\Maj_n \;|\;U) \gg \sig(\Maj_n \;|\;U)$ for any sequence of subsets with constant density.

\begin{quest}
Characterise sequences of Boolean functions  such that for any sequence of subsets $U_n$ with constant density $I(f_n \;|\;U_n) \gg \sig(f_n \;|\;U_n)$ holds, or where $I(f_n \;|\;U_n) \asymp \sig(f_n \;|\;U_n)$, respectively. 
\end{quest}

\subsection{Clue via distances between probability measures
}\label{cluedistmeas}

In this section, we introduce a somewhat different approach to measure the amount of information of a subset of coordinates about a Boolean function. 

Let us consider the usual configuration space $\left\{ -1,1 \right\}^V$ endowed with a probability measure $\mu$. Clearly, any Boolean function $ f: \left\{ -1,1 \right\}^V \longrightarrow \left\{ 0,1 \right\}$, after normalizing by $\E[f]$, can be interpreted as the density function of the  measure $\mu$ conditioned on the set of configurations $\{\omega \in    \left\{ -1,1 \right\}^V \; : \; f(\omega) =1 \}$. 

More generally, every $ f: \left\{ -1,1 \right\}^V \longrightarrow \R_{\geq 0}$ with $\E[ f]>0$ can be interpreted as a density, and can be used to define another probability measure on the same space, by
\beq \label{density}
\nu[\omega]: = \frac{1}{\E[ f]}f(\omega)\mu[\omega], \quad \omega \in \left\{ -1,1 \right\}^n.
\eeq

Now that we have identified our function with a probability measure, we can express clue in terms of distances of probability measures. We will consider three possible metrics: the total variation distance,  $L^2$ distance and in Section \ref{clueentropy} we use information theoretic distances.

We fix a (non-trivial) Boolean function $f$ and introduce the notation $\mu^1$ and  $\mu^0$ for the measures $\mu$ conditioned  on the set $\{f(\omega) =1 \}$ and $\{f(\omega) =0 \}$, respectively. Furthermore, let $ \mu[f(\omega) =1] = p$. So, we have $\mu = (1-p)\mu^0+ p\mu^1$ and, as a consequence,
\beq \label{Booledensity}
\frac{d\mu^1}{d\mu}  = \frac{1}{p}f \quad\text{and }\quad   \frac{d\mu^0}{d\mu} = \frac{1}{1-p}(1-f).
\eeq
We will need the measures $\mu|_U$, $\mu^1|_U$ and $\mu^0|_U$, which are the  marginals (projections) of the respective measures on the subset of coordinates  $U \subseteq V$. It is straightforward to check that
\beq \label{projdensity}
\frac{d\mu^1|_U}{d\mu|_U}  = \frac{1}{p} \E[f|\mathcal{F}_{U}]  \quad\text{and }\quad \frac{d\mu^0|_U}{d\mu|_U}  = \frac{1}{1-p}(1- \E[f|\mathcal{F}_{U}]).
\eeq
In addition, we still have $\mu|_U = (1-p)\mu^0|_U + p\mu^1|_U$. It is also worth noting that the first equalities in both~(\ref{Booledensity}) and~(\ref{projdensity}) work not only in the Boolean case, but for any nonnegative $f$, using~(\ref{density}) and $p = \E[f]$. Moreover, if  $f: \left\{ -1,1 \right\}^V \longrightarrow [0,1]$, then $1-f$ is again non-negative, and along with  $\mu^1$ one can define  $\mu^0$ as well, according to the formulas \eqref{Booledensity}.

We introduce two meaningful ways of measuring the clue, irrespective of the particular notion of distance $D(\cdot,\cdot)$ between probability measures.  In the first version, the total information content of the function is measured through the distance of $\mu^1$ from the original measure $\mu$:
\beq \label{asymclue}
\clue_{\text{asym}}^D (f\;|\;U) := \frac{D(\mu|_U, \mu^1|_U)}{D(\mu, \mu^1)}.
\eeq
In the alternative, symmetric version, which is the setup used in \cite[Chapter 16]{Pe},  we express the information content with the distance between the measures $\mu^1$ and  $\mu^0$:
\beq \label{symclue}
\clue_{\text{sym}}^D (f\;|\;U) := \frac{D(\mu^0|_U, \mu^1|_U)}{D(\mu^0, \mu^1)}.
\eeq
This version has the (desirable) general property that $\clue_{\text{sym}}^D(\1_{A}\;|\;U)= \clue_{\text{sym}}^D(\1_{A^{c}}\;|\;U)$ for any $U \subseteq V$. The asymmetric version has the advantage that it works not only for $f: \left\{ -1,1 \right\}^V \longrightarrow [0,1]$, as the symmetric one, but extends to all functions $ f: \left\{ -1,1 \right\}^n \longrightarrow \R_{\geq 0}$. 

In the sequel we shall discuss the $L^2$ and $L^1$ distances between probability measures. The squared $L^2$  distance between two measures $\nu, \theta \ll \mu$ is given by
\beq \label{D2}
D_{2}^2(\nu, \theta) :=  \frac{1}{4}\int_{\Omega}{(g-h)^2 d\mu},
\eeq
where $g = d\nu/d\mu$ and $h = d\theta/d\mu$. Accordingly, by \eqref{Booledensity}  we have  
$$
D_{2}^2(\mu, \mu^1) = \frac{1}{4}\int_{\Omega}{\left(1-\frac{f}{p}\right)^2 d\mu}= \frac{\Var(f)}{4p^2},
$$
where $p = \E[f]$. Similarly, \eqref{projdensity} gives that
$$
D_{2}^2(\mu|_U, \mu^1|_U) = \frac{1}{4}\int_{\Omega}{\left(1-\frac{ \E[f|\mathcal{F}_{U}]}{p}\right)^2 d\mu}= \frac{\Var( \E[f|\mathcal{F}_{U}])}{4p^2}.
$$

As for the symmetric version, again using \eqref{Booledensity}, we have 
\beq \label{symD2}
D_{2}^2(\mu^0, \mu^1) = \frac{1}{4}\int_{\Omega}{\left(\frac{f}{p} - \frac{1-f}{1-p}\right)^2 d\mu}= \frac{\Var(f)}{4p^2(1-p)^2}.
\eeq 
and, by \eqref{projdensity},
\beq \label{e.symD2cond}
D_{2}^2(\mu^0|_U, \mu^1|_U) = \frac{1}{4}\int_{\Omega}{\left(\frac{ \E[f|\mathcal{F}_{U}]}{p} - \frac{1- \E[f|\mathcal{F}_{U}]}{1-p}\right)^2 d\mu}= \frac{\Var( \E[f|\mathcal{F}_{U}])}{4p^2(1-p)^2}.
\eeq 
Thus we obtain that $\clue_{\text{sym}}^{L2}(f\;|\;U)=\clue_{\text{asym}}^{L2}(f\;|\;U) = \sqrt{\clue(f\;|\;U)}$, so it eventually does not matter whether we choose the symmetric or the asymmetric definition  of $\clue^{L2}$.

The total variation distance  (or $L^1$ distance)  is defined as
\beq  \label{TV}
D_{TV}(\nu, \theta) := \frac{1}{2}\sum_{\omega \in \Omega}{|\nu[\omega] - \theta[\omega]|} = \frac{1}{2}\int_{\Omega}{|g-h| d\mu}.
\eeq
Calculations similar to the $L^2$ distance case give that
$$
D_{TV}(\mu,  \mu^1) =  \frac{\E[|f-\E[f]|]}{2\E[f]},
$$
and 
$$
D_{TV}(\mu^0, \mu^1) =  \frac{\E[|f-\E[f]|]}{2\E[f]|\E[f]-1|}.
$$
So 
\beq \label{clueTV}
\clue^{TV}(f\;|\;U)= \frac{ \E[\left|\E[f|\mathcal{F}_{U}] - \E[f]\right| ]}{\E[|f-\E[f]|]},
\eeq
again irrespective of which of the two variants we use.

Without going into details we mention that one can again define a dual notion by $\sig^{TV}(f\;|\;U)= 1 -\clue^{TV}(f\;|\;U^{c}) $.

Let us compare $\clue$ and  $\clue^{TV}$ for Boolean functions.

\bpr \label{tvvsl2}
Let $ f: \left\{ -1,1 \right\}^V \longrightarrow \left\{ 0,1 \right\}$ with $\E[f] = p$  and  $U \subseteq V$. Then
$$
\frac{1}{L}\,\clue(f\;|\;U) \leq \clue^{TV}(f\;|\;U) \leq L\, \sqrt{\clue(f\;|\;U)}, 
$$
where $L = \max{\left\{\sqrt{\frac{1-p}{p}}, \sqrt{\frac{p}{1-p}} \right\}} $.
\epr

\bpf
We shall use $\clue_{\text{sym}}(f\;|\;U)$ here.
We have, by \eqref{Booledensity},
$$
 \left|\frac{d\mu^1}{d\mu}-\frac{d\mu^0}{d\mu}\right| = \left|\frac{f}{p}- \frac{1-f}{1-p}\right| = \frac{\left| f-p \right|}{p(1-p)} \leq \frac{1}{\min{\{p, 1-p \}}}.
$$
Since $f$ is Boolean, $\E[f|\mathcal{F}_{U}] \in [0,1]$, and in the same way, just using~\eqref{projdensity}, we have 
$$
\left|\frac{d\mu^1_U}{d\mu_U}-\frac{d\mu^0_U}{d\mu_U}\right| \leq \frac{1}{\min{\{p, 1-p \}}}.
$$

Therefore,
\begin{align*}
D^2_{2}(\mu^0, \mu^1) = 
\sum_{\omega} \mu[\omega]  \left|\frac{d\mu^1}{d\mu}(\omega)-\frac{d\mu^0}{d\mu}(\omega)\right|^2   
&\leq \frac{1}{\min{\{p, 1-p \}}} \sum_{\omega}{  \left|\frac{d\mu^1}{d\mu}(\omega)-\frac{d\mu^0}{d\mu}(\omega)\right| }\\
&=\frac{1}{\min{\{p, 1-p \}}}D_{TV}(\mu^0, \mu^1).
\end{align*}
In light of \eqref{projdensity}, it is clear that the same inequality holds for the pair $D_2^2(\mu^0|_U, \mu^1|_U)$ and  $D_{TV}(\mu^0|_U, \mu^1|_U)$. 

For an opposite bound, the Cauchy-Schwarz inequality gives
\begin{align*}
D_{TV}(\mu^0, \mu^1) &=\frac{1}{2}\sum_{\omega}{ \sqrt{\mu[\omega]}\sqrt{\mu[\omega]}
 \left|\frac{d\mu^1}{d\mu}(\omega)-\frac{d\mu^0}{d\mu}(\omega)\right| }\\
 &\leq\frac{1}{2} \sqrt{\sum_{\omega}{\mu[\omega]  \left|\frac{d\mu^1}{d\mu}(\omega)-\frac{d\mu^0}{d\mu}(\omega)\right|^2 }}= \frac12 D_{2}(\mu^0, \mu^1).    
\end{align*}

Since $f$ takes on only $0$ and $1$, we have $\Var(f) =  p(1-p)$, and thus~\eqref{symD2} tells us that  $D_2(\mu^0, \mu^1) =1/(2\sqrt{p(1-p)})$ and $2 \min\{p, 1-p \} D_2(\mu^0, \mu^1) = \min\left\{ \sqrt{\frac{p}{1-p}},\sqrt{\frac{1-p}{p}} \right\}= \frac{1}{L}$. Using this,
\begin{align*}
\frac{1}{L} \, \clue(f\;|\;U)&= \frac{ 2 \min\{p, 1-p \} D_2^2(\mu^0|_U, \mu^1|_U)}{D_2(\mu^0, \mu^1)}\\
& \leq  \clue^{TV}(f\;|\;U) \\
& \leq  \frac{D_2(\mu^0|_U, \mu^1|_U) }{2 \min \{p, 1-p \} D_2^2(\mu^0, \mu^1)}= L\,\sqrt{\clue(f\;|\;U)},
\end{align*} 
as claimed.
\epf

A sequence of Boolean functions  $f_n: (\Omega^{V_n},\,\p) \lora \left\{ -1,1 \right\}$ is called degenerate if $\Var(f_n) \to 0$, non-degenerate otherwise.
So Theorem \ref{cluethm} implies that, for a non-degenerate sequences of Boolean functions sparse reconstruction with respect to $\clue$ is equivalent to sparse reconstruction with respect to $\clue^{TV}$. In particular,  there is still no sparse reconstruction on product measures if we replace $\clue$  by $\clue^{TV}$.

\subsection{Clue via entropy} \label{clueentropy}


Our setup remains the same, but we formulate it in a somewhat different way. Let $\left\{ X_v : v \in V \right\}$ be a set of real-valued discrete random variables defined on a common probability space. As in the introduction, for $S \subseteq V$ we let $X_S  := \left\{ X_j : j \in S \right\}$.
 
The variables $\{X_v : v \in V\}$  play the role of the coordinates of Section~\ref{s.L2}. Indeed, alternatively we could talk about a measure $\mu$ on $\R^V$, and then $\p[ X_V \in B] = \mu[B]$ for a Borel-set $B \subseteq \R^V$. Again we consider a function  $f: \R^V \lora \R$ and let $Z= f(X_V)$. Our main interest is still the special case where the variables $X_v$ and $Z$ are binary valued variables (spins), 
but the arguments we present here work in this slightly more general framework.

In this section we are going to discuss an alternative way of measuring the amount of information a subset  $S \subseteq V$ of coordinates contains about the function $f$. We will use concepts from information theory and define an information-theoretic clue accordingly.

For the sake of completeness we start with some classical definitions. For a (possibly vector valued) random variable (or a probability distribution) entropy quantifies the amount of information in a sample.

\bde[Entropy]
Let $X$ be a discrete random variable. Then the entropy of $X$ is
$$
H(X) = -\sum_{x \in \mathrm{Ran}(X)}{\p[X=x] \log{\p[X=x] } }. 
$$
\ede

We will also need the concept of conditional entropy. The entropy of $X$ conditioned on the random variable $Y$ expresses how much randomness remains in $X$ on average if we learn the value of $Y$.

\bde[Conditional entropy]
Let $X$ and $Y$ be two discrete random variables defined on the same probability space. The conditional entropy of $X$ given $Y$, 
is
$$
H(X \given Y) = \sum_{y \in \mathrm{Ran}(Y)}  {\p[Y = y]} H(X  \given Y = y),
$$
where $H(X  \given Y = y)$ is the entropy of the distribution of $X$ conditioned on the event $\{Y=y\}$.
\ede

The mutual information quantifies the common information present in two variables. In a way it measures how far the joint distribution of the two variables is from being independent.

\bde [Mutual Information] \label{muti}
Let $X$ and $Y$ be two discrete random variables defined on the same probability space. Suppose that $H(X)$ and $H(Y)$ are both finite. Then the mutual information between $X$ and $Y$ is:
\beq
I(X : Y ) = H(X) + H(Y) - H(X,Y ) = H(X) - H(X \given Y).
\eeq
\ede

Now comes the definition of clue in this framework. 

\bde[I-clue] \label{Infcluedef}
Let $\{X_v :\;v \in V\}$ be a finite family of discrete real valued random variables defined on the same probability space, and for some  $f: \R^V \lora \R$ let us consider the random variable $Z = f(X_V)$. The information theoretic clue (I-clue) of $f$ with respect to $U \subseteq V$  is
$$
\clue^{I}(f\;|\;U)= \frac{I(Z : X_U )}{I(Z : X_V )} = \frac{I(Z : X_U )}{H(Z)}.
$$
\ede

Note that if $Z$ is $X_U$-measurable then $H(Z \given X_U)=0$, and therefore $I(Z:X_U ) = H(Z)$, while if $Z$ is independent from  $X_U$ then $I(Z: X_U)=0$, in accordance with what we expect from a clue-type notion.

As for the cases discussed before, here too we can introduce the dual (which expresses again how much information we are missing if we do not know the coordinates in $U$):
$$
\sig^{I}(f\;|\;U)= 1 - \frac{I(Z : X_{U^{c}} )}{H(Z)} = \frac{H(Z \;|\; X_{U^{c}} )}{H(Z)} .
$$
In the sequel we show that $\clue^{I}$ can  also be interpreted via distances of probability measures. For this, we will need the following definition.

In information theory the following concept is used to measure distance between probability measure. 

\bde[Kullback-Leibler divergence]
Let $\mu$ and $\nu$  be probability measures on the same discrete probability space $\Omega$, where $\nu \ll \mu$. The relative entropy between  $\nu$ and $\mu$ is
$$
D_{KL}(\nu \,\|\, \mu) = -\sum_{x \in \Omega}{\nu(x) \log{\frac{\mu(x)}{\nu(x)} } }. 
$$
\ede

The KL divergence, although it also means to express a concept of distance between two distributions, is not a metric. In particular, it is neither symmetric, nor does it satisfy the triangle inequality.

We can introduce yet another notion of clue for any nonnegative function $f$ with $\E[f \log f] < \infty$. As in Section~\ref{cluedistmeas}, we now interpret $f/\E[f]$ as the density  (or Radon-Nikodym derivative) of some measure $\nu$ .

Let  $ \mathsf{Ent}(f) : = \E[f \log{f}] - \E[f]\log\E[f]$; the expectation is taken with respect to $\mu$, and in case  $f(\omega) =0$, we have, by continuity, $f(\omega) \log f(\omega) = 0$. One can easily verify that the following holds: 
$$
D_{KL}(\nu\,\|\,\mu ) = \mathsf{Ent}\left(\frac{d\nu}{d\mu}\right)  = \mathsf{Ent}\left( \frac{f}{\E[f]}\right) = \frac{1}{\E[f]} \mathsf{Ent}( f).
$$
We can use this to introduce an  asymmetric notion of clue (in the sense of  Section \ref{cluedistmeas}):
$$
\clue^{KL}(f\;|\;U):= \frac{D_{KL}(\nu_{U} \,\|\,\mu_{U})}{D_{KL}(\nu \,\|\,\mu)} = \frac{\mathsf{Ent}(\E[f \given U])}{\mathsf{Ent}(f)}.
$$
Observe that there is no symmetric version of clue here. The problem is that in the Boolean case if $\mu^0$ and $\mu^1$ are the projection of $\mu$ on $\{ f = 1 \}$ and  $\{ f = 0 \}$, respectively as in Section \ref{cluedistmeas}, then $D_{KL}(\mu^0 \,\|\ \mu^1) = \infty$, since  $\mu^0$ and $\mu^1$ are singular.

Note the similarity of $\clue^{KL}$ as expressed in terms of $\mathsf{Ent}$ and the original, $L^2$ version of $\clue$ expressed in terms of variance. Indeed, $\mathsf{Ent}(f)$ and $\Var(f)$ together with the respective concepts of clue can be examined in the general framework of $\Phi$-entropies  (see, for example, \cite[Chapters 14 and 15]{BLM}). The main idea is that for any convex function $\Phi: \R \lora \R $ one can assign a respective $\Phi$-entropy for every integrable random variable $X$, by
$$
H_{\Phi}(X) = \E[\Phi(X) ] - \Phi(\E[X]).
$$
It turns out that under some general analytic conditions on $\Phi$, many important properties we require from an information measure remain valid for $H_{\Phi}(X)$ (for example, it is always non-negative because of Jensen's inequality). In particular, we get $H_{\Phi}(X) =\Var(X)$ when $\Phi(x) = x^2$ and  $H_{\Phi}(X) = \mathsf{Ent}(X)$ with $\Phi(x) = x \log x$. 

We note that $\clue^{I}$ can also be interpreted via distances between probability measures. A standard calculation shows that the mutual information can be expressed through Kullback-Leibler (KL) divergence as follows:
\beq\label{e.IKL}
I(Z : X_{[n]} ) = \E_{Z}\big[D_{KL} \big(\mu|Z \,\big\|\,  \mu\big) \big]
\eeq
Here $\mu|Z$ denotes the measure $\mu$ conditioned on possible values of $Z = f(\omega)$. In particular,  when \eqref{Booledensity} and thus \eqref{projdensity} are satisfied (for example, when $f$ is Boolean), one can introduce another information theoretic distance (defined in \cite[Chapter 16]{Pe} for $p =1/2$) as follows: 
\beq\label{e.DKL}
D_{I}^{ \mu}(\mu^1,\mu^0 ) :=\E_{Z}\big[D_{KL} \big(\mu|Z \,\big\|\,  \mu\big) \big] = p\mathsf{Ent}\left( \frac{ d\mu ^1}{d\mu}\right) + (1-p) \mathsf{Ent}\left( \frac{ d\mu ^0}{d\mu} \right) .
\eeq
Thus, using the distance $D_{I}^{ \mu}$, by \eqref{e.IKL} the $\clue^{I}$ can be written as follows:
$$
\clue^{I}(f\;|\;U)=\frac{\E\big[D_{KL} \big(\mu |_U |Z \,\big\|\,  \mu|_U\big)\big]}{\E_{Z}\big[D_{KL} \big(\mu|Z \,\big\|\,  \mu\big) \big]} 
= \frac{D_{I}^{ \mu|_U}(\mu^0|_U, \mu^1|_U)}{D_{I}^{ \mu}(\mu^0, \mu^1)}.
$$


The following proposition shows that for non-degenerate sequences  of Boolean functions,  sparse reconstruction  with respect to $\clue$ and $ \clue^{I}$ are equivalent.

\bpr \label{clueequ}
Let $\mu$ be a measure on $ \{-1,1\}^{V}$ and $\{\omega_v : v \in V\}$ a spin system distributed according to $\mu$.
Let  $f: \{-1,1\}^{V} \longrightarrow \{0,1\}$ and  $U \subseteq V$. Then,
\beq
\frac{1}{2}\E[f]^2 \, \E[1-f]^2 \, \clue^{I}(f\;|\;U) \leq \clue(f\;|\;U) \leq\frac{4 \log{2}}{p_{\mathrm{min}}} \clue^{I}(f\;|\;U),
\eeq
where $p_{\mathrm{min}} := \min(\mu[f=1], \mu[f=0])$.
\epr

\bpf
We will bound mutual information by conditional variance, and vice versa. For the first of these, we use the setup of Section~\ref{cluedistmeas} and follow the idea sketched in the proof of \cite[Lemma 16.5 (ii)]{Pe}. Recall that $\mu^1$ and  $\mu^0$ are the measures $\mu$ conditioned  on the sets $\{f(\omega) =1 \}$ and $\{f(\omega) =0 \}$, respectively.  We are going to show that
$$
D_I(\mu ^0, \mu ^1)\leq  4 D_2^2(\mu ^0, \mu ^1).
$$
We introduce the shorthand notation $g := \frac{ d\mu ^1}{d\mu} = \frac{ f}{p} $ and $h := \frac{ d\mu ^0}{d\mu} = \frac{1-f}{1-p}$. Introducing  $\psi :=  pg - (1-p)h$, we have  $\frac{1+\psi}{2} = pg$ and $\frac{1-\psi}{2} = (1-p)h$.
This allows us to write 
$$
D_I(\mu ^0, \mu ^1) = p\mathsf{Ent}(g) + (1-p)\mathsf{Ent}(h) = \int_{\Omega}{\frac{1+\psi}{2} \log{\frac{1+\psi}{2p}} + \frac{1-\psi}{2} \log{\frac{1-\psi}{2(1-p)}}d\mu},
$$
where $\log$ denotes the natural logarithm. Now we divide the integral as follows:
\begin{align*}
   D_I(\mu ^0, \mu ^1)&=\int_{\Omega}{\frac{1+\psi}{2} \log{(1+\psi)} + \frac{1-\psi}{2} \log{(1-\psi)} \, d\mu} \\
   & \qquad - \int_{\Omega}{p\log{(2p)} g  + (1-p)\log{(2(1-p))} h \, d\mu} .
\end{align*}
Using that $\int_{\Omega}{g d\mu}  = \int_{\Omega}{h d\mu}=1$ and  that the entropy of a binary-valued random variable can be at most $\log 2$, we conclude that the second term is always negative: 
\begin{align*}
\int_{\Omega}{p \log({2p}) g + (1-p) \log({2(1-p)})h \, d\mu} 
&=  \log 2 + p\log{p}+ (1-p)\log{(1-p)}\\
& = \log 2 - H(Z)\geq 0. 
\end{align*}
In case $p \ge 1/2$, we have $|\psi| = p\left |g - \frac{1-p}{p} h\right | \leq p |g-h|$, because the supports of $g$ and $h$ are disjoint. Similarly, for $p \leq 1/2$ we have $|\psi| \leq (1-p)|g-h|$.  Altogether, $|\psi| \leq \max(p, 1-p)|g-h|$.  Moreover, using $\log(1+x)\leq x$   we obtain that
\begin{align*}
D_I(\mu ^0, \mu ^1) &\leq \int_{\Omega}{\frac{1+\psi}{2} \log{(1+\psi)} + \frac{1-\psi}{2} \log{(1-\psi)}\,d\mu} \leq \int_{\Omega}{\psi^2 d\mu} \\
&\leq \max(p, 1-p)^2 \int_{\Omega}{(g-h)^2 d\mu} \leq 4 D_2^2(\mu ^0, \mu ^1) = \frac{\Var(f)}{p^2(1-p)^2}.
\end{align*}

Now we apply this inequality to the measures $\mu|_U=p\mu^1|_U+(1-p)\mu^0|_U$ defined in~(\ref{projdensity}), using $\E[Z \,|\, \mathcal{F}_U]$ with $Z = f(\omega)$.  By~\eqref{e.IKL} and~\eqref{e.DKL}, the left hand side becomes $D_I(\mu^0|_U, \mu ^1|_U)  = \E\big[D_{KL} \big(\mu |_U |Z \,\big\|\,  \mu|_U\big)\big]= I(Z : \omega_U )$. On the right hand side, $D_{2}^2(\mu^0|_U, \mu^1|_U)$ can be rewritten using~\eqref{e.symD2cond}. Thus, we get
\beq \label{I<Var}
I(Z : \omega_U ) \leq  \frac{\Var(\E[Z \,|\, \mathcal{F}_U])}{\E[Z]^2(1-\E[Z])^2}.
\eeq

We turn to the second comparison. We show that, more generally, whenever $f: \{-1,1\}^{n} \longrightarrow [-1,1]$, we have
\beq \label{var<I}
\Var(\E[Z \,|\, \mathcal{F}_U]) \leq 2 I(Z : \omega_U ).
\eeq
Our argument follows \cite[Lemma 4.4]{Ta}. 
First we fix some notations. Let $z$ be in the range of $f$ and $u \in  \{-1,1\}^{U} $. Then  
$$
p_{z} := \mu[Z = z], \quad p_{u} := \mu[ \omega_U = u], \quad p_{z|u} := \mu[Z = z \,|\, \omega_U = u ].
$$
Now, with this notation we have 
$$
\Var(\E[Z \,|\, \mathcal{F}_U]) = \sum_{u \in   \{-1,1\}^{U}}{p_u (\E[Z] -\E[Z \,|\,  \omega_U = u] )^2},
$$
and for a fixed $u \in  \{-1,1\}^{U} $
$$
\big(\E[Z] -\E[Z \,|\,  \omega_U = u] \big)^2 = \sum_{z }{(p_{z}z -p_{z|u}z )^2} \leq \sum_{z}{(p_{z} -p_{z|u} )^2}.  
$$
So we get that 
\beq \label{varestimate}
\Var(\E[Z \,|\, \mathcal{F}_U])  \leq \sum_{u \in   \{-1,1\}^{U}}p_u{\sum_{z}{(p_{z} -p_{z|u} )^2}}.
\eeq
With the notation $h(x):= -x \log x$ for $x \in [0,1]$ (where $h(0) := 0$) we can write the mutual information as  
\beq \label{miexpansion}
I(Z : \omega_U ) = H(Z) - H(Z|\omega_U) = \sum_{z}{\left(h(p_{z}) - \sum_{u \in   \{-1,1\}^{U}}{ p_u h(p_{z|u})}\right)}.
\eeq
Using  linear Taylor expansion for $h(p_{z|u})$ around $p_z$ with error term, we get the following estimate:
$$
h(p_{z|u}) =  h(p_{z}) + h'(p_{z})(p_{z|u} -p_{z}) -\frac{1}{2p^{*}_{z|u}}(p_{z|u} -p_{z})^2,
$$
with some $p^{*}_{z|u}$ between $p_{z|u}$ and $p_{z}$, using for the error term that  $ h''(x) = -\frac{1}{x}$.
Substituting this estimate into \eqref{miexpansion},  we observe that the terms with $h'(p_{z})$ cancel, since for any $z \in \{0,1\}$ we have  $\sum_{u \in   \{-1,1\}^{U}}{ p_u (p_{z|u} - p_z)} = p_z-p_z = 0$.  Therefore we obtain that
$$
\sum_{u \in  \{-1,1\}^{U}}{p_u\sum_{z}{\frac{(p_{z} -p_{z|u} )^2}{p^{*}_{z|u}}}} = 2 I(Z : \omega_U ).
$$
As $0<p^{*}_{z|u} <1$ we can conclude, using \eqref{varestimate} that
$$
\Var(\E[Z \,|\, \mathcal{F}_U])  \leq \sum_{u \in  \{-1,1\}^{U}}{p_u\sum_{z}{\frac{(p_{z} -p_{z|u} )^2}{p^{*}_{z|u}}}} \leq 2  I(Z : \omega_U ).
$$

Finally, in order to get a stronger bound we shall use a sharper inequality between the entropy and the variance in the denominators. Observe that for Boolean functions we have $\Var(Z) = p(1-p) \geq \frac{p_{\mathrm{min}}}{2}$. At the same time (again because $f$ is Boolean) $H(Z) \leq \log 2$. Thus we have
\beq \label{entvsvar}
H(Z) \leq  \log 2 \leq  \frac{2\log 2 }{p_{\mathrm{min}}}\Var(Z).
\eeq

$$
$$


Now we get the first inequality of the statement by using  \eqref{I<Var} for the numerator and \eqref{var<I} with  $U=V$ for the denominator:
$$
\frac{1}{2} p^2(1-p)^2 \clue^{I}(f\;|\;U) = \frac{p^2(1-p)^2 I(Z : \omega_U )}{2 I(Z : \omega_V )} \leq \frac{\Var(\E[Z \,|\, \mathcal{F}_U])}{\Var(\E[Z \,|\, \mathcal{F}_V])} = \clue(f\;|\;U).
$$ 
For the second inequality we apply \eqref{var<I} for the numerator and \eqref{entvsvar} for the denominator:
\[
\clue(f\;|\;U) = \frac{\Var(\E[Z \,|\, \mathcal{F}_U])}{\Var(Z)}  \leq \frac{\log{4} }{p_{\mathrm{min}}} \frac{2 I(Z : \omega_U )}{ H(Z )} =  \frac{4\log{2} }{p_{\mathrm{min}}} \clue^{I}(f\;|\;U).
\qedhere
\]

\epf

\section{Sparse reconstruction with respect to I-clue and KL-clue}\label{s.infoclue}

In this section we show some analogues of Theorem \ref{cluethm2} for the I-clue and KL-clue. We note that the following theorem,  as well as the definition of I-clue only works well in the discrete case, as the continuous counterpart of entropy, differential entropy has some drawbacks (for example, it can be negative). The notation and setup follows Subsection \ref{clueentropy}.

\bth \label{entrcluethm}
Let $\left\{ X_v  : v \in V \right\}$ be discrete valued iid~random variables with finite entropy. Let $f: \Omega^V \longrightarrow \R$ be a transitive function and $Z = f(\left\{ X_v  : v \in V \right\})$. Then
\beq
\clue^{I}(f\;|\;U)\leq \frac{|U|}{|V|}.
\eeq
\eth

For the proof we will use the following well-known inequality which finds numerous applications in combinatorics. For a proof see, for example, \cite[Theorem 6.28]{LP}.
 
\bth[Shearer's inequality \cite{Shearer}] \label{shearer}
Let $X_1, X_2, \dots X_n$ random variables defined on the same probability space. Let $S_1,S_2, \dots, S_L$ be subsets of $[n]$ such that for every $i \in [n]$ there are at least $k$ among $S_1,S_2, \dots, S_L$ containing $i$. Then
$$
k H(X_{[n]}) \leq \sum_{j=1}^{L}{H(X_{S_j})}.
$$
\eth

First we need the following consequence of Shearer's inequality.

\bl \label{InfoShearer} Suppose $X_1, X_2, \dots, X_n$ are independent, and $Z=f(X_1,\dots,X_n)$.
Let $S_1,\dots, S_L$ be a system of subsets of $ [n]$ such that each $i \in  [n]$ appears in at most $k$ sets. Then
\beq
\sum_{j=1}^{L}{I(Z:X_{S_j} )}\leq k I(Z: X_{[n]}). 
\eeq
\el

\bpf
Without loss of generality we can assume that each $i$ appears in exactly $k$ sets. Indeed, if this is not the case, we can always add some additional subsets so that this condition is satisfied. While adding new sets the right hand side of the inequality does not change and the left hand side can only increase.

Since the  variables $X_i$ are independent: 
\beq \label{indep}
\sum_{j=1}^{L}{H(X_{S_j})} = \sum_{j} {\sum_{i \in S_j}{H(X_i)}} = k {\sum_{i \in [n]}{H(X_i)}}=k H(X_{[n]}).
\eeq 
On the other hand, using Shearer's inequality,
\beq \label{searer}
- \sum_{j=1}^{L}{H(X_{S_j}|Z)} \leq -k H(X_{[n]}|Z) .
\eeq
Using that $I(Z: X_{S_j} )= H(X_{S_j}) - H(X_{S_j}|Z)$, adding up \eqref{indep} and \eqref{searer} completes the proof. 
\epf

Now the proof of the clue-theorem:

\bpf
Recall that we have the iid measure $\p$, and the function $f$ is invariant under some transitive action of a group $G$. Let $U \subseteq V$ be arbitrary. Then, for each $g \in G$, 
$$
I(Z : X_U) = I(Z : X_{U^{g}})\,,
$$
where $U^{g} = \left\{ ug\; :\;g \in G \right\}$.

Observe that $v \in U^{g} \;\Longleftrightarrow\; \exists\, u \in U \textrm{ with } vg^{-1}=u$. For each pair of $v \in V$ and $u \in U$ there are $|G_{v}|$ such $g$, where $G_{v}$ is the stabilizer subgroup of $G$ at $v$. (Since the action is transitive such a $g$ exists, moreover the cardinality of the stabilizer subgroup $G_{v}$ is the same for every $v \in V$.) The conclusion is that each $v \in V$ appears in exactly $|U| \cdot |G_{v}|$ translated versions of $U$. 
Applying Lemma \ref{InfoShearer} gives 
$$
|G| \, I(Z: X_U)=\sum_{g \in G}{I(Z : X_{U^{g}} )}\leq |U| \, |G_{v}| \, I(Z : X_{V}) = |U| \, |G_{v}|\, H(Z),
$$
which is what we wanted since $|G| = n|G_{v}|$ by the orbit-stabilizer theorem.
\epf

Observe that for any  non-degenerate sequence of Boolean functions,  sparse reconstruction with respect to I-clue is equivalent to sparse reconstruction with respect to the original, $L^2$ version (irrespective of the underlying measure). This follows from Proposition \ref{clueequ}. Nevertheless, in case $\{f_n\}$ is degenerate Boolean, or non-Boolean, Proposition \ref{clueequ} does not help us compare the sequences of clues and  I-clues. This raises the following question:

\begin{quest}
Is there a sequence  of functions $f_n \in L^2( \{-1,1\}^{V_n},\pi_n^{\otimes V_n})$  and a corresponding sequence of subsets $U_n \subseteq V_n$ such that 
\begin{enumerate}
    \item $\clue^{I}(f_n\;|\;U_n) \ll \clue(f_n\;|\;U_n)$
    \item $\clue^{I}(f_n\;|\;U_n) \gg \clue(f_n\;|\;U_n)$?
\end{enumerate}
What is the answer if we ask $f_n$ to be Boolean for all $n \in \N$?
What is the answer if we allow for non-product measures on $\{-1,1\}^{V_n}$?
\end{quest}

It is remarkable that for product measures we have the same inequality for the clue and I-clue of general (possibly degenerate or $\R$-valued) sequences of  transitive functions. In particular, we emphasize that Theorem \ref{cluethm2} and Theorem \ref{entrcluethm} do not imply one another.



\medskip
Interestingly enough, along the same logic one can prove the respective version of Theorem \ref{cluethm2} and Theorem \ref{entrcluethm} for $\clue^{KL}$.  We should emphasize that, in contrast with mutual information, relative entropy is a concept that remains meaningful for continuous random variables as well. So Theorem \ref{KLcluethm} holds for all product measures, just like Theorem~\ref{cluethm2}. The proof relies on the following Shearer-type inequality:

\bl \label{InfoShearer1}
Let $\p$ be a product measure, and $\mu$ another probability measure on the same space satisfying $\mu \ll \p$. 
Let $S_1,\dots, S_L$ be a system of subsets of $ V$ such that each $i \in  V$ appears in at most $k$ sets. Then
$$
\sum_{j=1}^{L} D(\mu_{S_i} \,\|\,\p_{S_i}) \leq k D(\mu \,\|\,\p ).
$$
\el

In our application, of course $\mu$ is the measure with density $f$. It is easy to recognise that Lemma \ref{InfoShearer1} is a close relative of Lemma \ref{InfoShearer}. The proof of this lemma is also a straightforward consequence of Shearer's inequality (Theorem \ref{shearer}); for a proof see \cite{Ga}. The corresponding clue theorem follows in the same way as Lemma \ref{InfoShearer} implies Theorem \ref{entrcluethm}.

\bth \label{KLcluethm}
Let $\left\{ X_v  : v \in V \right\}$ be $\Omega$-valued iid random variables. Let $f: \Omega^n \longrightarrow \R$ be a transitive function and $Z = f(\left\{ X_v  : v \in V \right\})$. Then
\beq
\clue^{KL}(f\;|\;U)\leq \frac{|U|}{n}.
\eeq
\eth

It is worth noting that for sequences of transitive Boolean functions on the hypercube there is no sparse reconstruction, no matter which version of clue we wish to choose. Indeed, $W$ (witness) is dominated by $\clue$  (see Section \ref{siginf}), so  Theorem~\ref{cluethm} applies. As for $\clue^{TV}$, Proposition~\ref{tvvsl2} implies that whenever $\clue$ converges to $0$,  so does  $\clue^{TV}$, and again we can use Theorem~\ref{cluethm}. For $\clue^{I}$ and $\clue^{KL}$, this has been shown in the present section (Theorems~\ref{entrcluethm} and~\ref{KLcluethm}).

\section{Sparse reconstruction and cooperative game theory}\label{s.coop}

The field of cooperative game theory (for an introduction see, for example, \cite{BDT} or \cite{PS}) starts with the following setup: there is a set of players which we denote by $V$ here (to be consistent with our previous notation), and the game is defined by assigning a nonnegative real number $v(S)$ to every subset $S$ of the players. Usually it is assumed that $v(\emptyset)=0$ and that $v(S) \leq v(T)$, whenever $S \subseteq T$. The function $v : 2^{V} \longrightarrow \R$ is referred to as the characteristic function. This aims to model a situation where individuals can gain profit, and the profit may increase in case a group of individuals cooperates by forming a coalition. Thus  $v(S)$ is the joint payoff of the individuals in $S$ provided that they cooperate.

Cooperative game theory is  concerned with finding a fair distribution of the payoff given the characteristic function $v$.  One of these concepts is the Shapley value, introduced in \cite{Sh}, which distributes the payoff based on the average marginal contribution of the individuals.

\bde [Shapley value]
\beq
\phi_i(v) = \frac{1}{|V|}\sum_{S \subseteq V\setminus\{ i \}}{\frac{v(S\cup \{ i \})-v(S)}{{{|V|-1} \choose |S|}}}.
\eeq
\ede

A straightforward calculation shows that $\sum_{i \in V}{\phi_i(v)} = v(V)$. So the Shapley value is indeed the distribution of the payoff of the grand coalition. In general, this is not true for smaller coalitions $S \subset V$.

Observe that  for a given $f: \{-1,1\}^{V} \longrightarrow \{-1,1\}$ we can define a cooperative game via  $ v_{f}(U): = \Var[\E[f \,|\, \mathcal{F}_U]]$ for any $U \subseteq V$. Besides fitting the mathematical definition, it also fits into the interpretation of the theory. It is a sort of an information game, where the payoff depends on how accurately we know a piece of information (represented by the value of the function). Each individual possesses one piece of information (the value of the corresponding coordinate) but only together they determine the valuable piece of information.

In the proof of Theorem \ref{cluethm} we introduced the random element $X$ of the index set, which is a uniformly random element of the Spectral Sample. In fact, $X$ is distributed according the Shapley value.

\bpr \label{Xisshapley}
Let $f: \{-1,1\}^{V} \longrightarrow \R$. Then 
$$
\frac{\phi_i( v_{f})}{v_{f}(V)} = \p[X = i].
$$
\epr

\bpf
Without loss of generality we may assume that $\Var(f) =1$. Let $n = |V|$. First, observe that
$$
\p[X=u] = \sum_{S \ni u}{\widehat{f}(S)^2\frac{1}{|S|}}.
$$
Now we calculate $\phi_i( v_{f})$ via Fourier-Walsh expansion and show that it equals $\p[X=u]$. Using that $v_{f}(S) = \sum_{T \subseteq S}{\widehat{f}(T)^2}$ we get that
$$
\phi_i(v) = \frac{1}{n}\sum_{S \subseteq V \setminus\{ i \}}{\frac{\sum_ {T \subseteq S} {\widehat{f}(T\cup \{ i \})^2}}{{{n-1} \choose |S|}}} =  \frac{1}{n} \sum_{T \subseteq V\setminus\{ i \}}{\widehat{f}(T\cup \{ i \})^2 \sum_{S \subseteq [n]\setminus\{ i \} :\\ T \subseteq S} {\frac{1}{{{n-1} \choose |S|}}} } .
$$
For a fixed $T$ there are  ${{n-1 - |T|} \choose k-|T|}$ $k$-element subsets $S$ which contain $T$. Therefore we have 
$$
\phi_i(v) = \frac{1}{n} \sum_{T \subseteq V\setminus\{ i \}}{\widehat{f}(T\cup \{ i \})^2 \sum_{k= |T|}^{n-1} {\frac{{{n-1 - |T|} \choose k-|T|}}{{{n-1} \choose k}}} } 
$$
With some elementary manipulation of the binomial coefficients we get that
$$
\frac{{{n-1 - |T|} \choose k-|T|}}{{{n-1} \choose k}} = \frac{{{k} \choose |T|}}{{{n-1} \choose |T|}}.
$$
Now we apply the so called hockey-stick identity --- $\sum_{k= |T|}^{n-1} {{{k} \choose |T|} } = {{n} \choose |T|+1}$ --- and we get the desired formula.
\epf

Given how naturally the Shapley value arises in the proof of Theorem \ref{cluethm}, it is perhaps not surprising that there is a proof that does not use Fourier-Walsh transform, only simple concepts from cooperative game theory and combinatorics. The advantage of this approach is that it makes  transparent the reasons behind the striking similarities between Theorem~\ref{cluethm2} and Theorem~\ref{entrcluethm}. It should also be noted that this approach entails both theorems.

We need to introduce another concept of fair distribution which is related to our topic. The core, introduced in \cite{Gi}, defines those distributions of the profit in which every coalition of players gets in total at least as much as they deserve (according to the characteristic function).

\bde[Core]
The core of a cooperative game $v$ with set of players $V$ is the set  $C(v) \subseteq  \R^{|V|}$ such that $x \in C(v)$ if and only if
$$
\sum_{i \in V}{x_i} = v(V),
$$
and, for every $S \subset V$,
$$
\sum_{i \in S}{x_i} \geq v(S).
$$
\ede

We have the following simple observation.

\bl \label{coreshap}
Let $v$ be a game on the finite set $V$, and let $\Gamma$ be a group acting on $V$ such that $v$ is invariant under the action of $\Gamma$. If the Shapley value vector $\phi(v)$ is in the core $C(v)$, then, for every $S \subseteq V$,
$$
v(S) \leq  \frac{|S|}{\min_{ i \in V}{|\Gamma \cdot i|}}v(V).
$$
\el

\bpf
On the one hand, for every $j \in \Gamma \cdot i$ (recall that $ \Gamma \cdot i$ is the $\Gamma$--orbit of $i$), we have $\phi_i(v) = \phi_j(v)$, by symmetry. So we have
$$
|\Gamma \cdot i|\phi_i(v) = \sum_{j \in \Gamma \cdot i}{\phi_j(v)} \leq \sum_{j \in V}{\phi_j(v)} =   v(V).
$$
That is, 
$$
 \phi_i(v)\leq \frac{v(V)}{|\Gamma \cdot i|}.
$$
Using that  $\phi(v) \in C(v)$, we get that
$$
v(S) \leq  \sum_{i \in S}{\phi_i(v)} \leq  \frac{|S|}{\min_{ i \in V}{|\Gamma \cdot i|}}v(V),
$$
as claimed.
\epf

We are going to show that a class of cooperative games, the so-called convex games, satisfy the conditions of Lemma \ref{coreshap}. This concept was also first studied by Shapley (see \cite{Sh71}).

\bde [Convex games]
A cooperative game $v$ is convex if the characteristic function is supermodular. That is, for every subset of players $S, T \subseteq [n]$,
\beq \label{supermodular}
v(S)+v(T) \leq v(S \cup T)+v(S \cap T).
\eeq
\ede

The subgame $v_U$ denotes the game $v$ with its domain restricted to the subset  $U \subseteq [n]$.

\bl \label{supmod}
Let $v$ be a supermodular cooperative game. Then the  vector $\phi(v)$ is in the core $C(v)$.
\el

\bpf
We first show that $S\subseteq T$ implies that $\phi_i(v_S)  \leq \phi_i(v_T)$. It  suffices to prove this when $T= S\cup\{ j \}$. Let $|S|= k$.
 
It is a straightforward calculation to verify that, for any $l < k$,
$$
 \frac{1}{k}\frac{1}{{k-1 \choose l}} =  \frac{1}{k+1}\left(\frac{1}{{k \choose l}} +\frac{1}{{k \choose l+1}} \right),
$$
and therefore one can write
\begin{align*}
 \phi_i(v_S)&=  \frac{1}{k}\sum_{L \subseteq S\setminus\{ i \}}{\frac{v(L\cup \{ i \})-v(L)}{{{k-1} \choose |L|}}}\\ 
 &= \frac{1}{k+1}\left(\sum_{L \subseteq S\setminus\{ i \}}{\frac{v(L\cup \{ i \})-v(L)}{{{k} \choose |L|}}} + \sum_{L \subseteq S\setminus\{ i \}}{\frac{v(L\cup \{ i \})-v(L)}{{{k} \choose |L|+1}}} \right).
\end{align*}
At the same time,
$$
 \phi_i(v_T) = \frac{1}{k+1} \left( \sum_{L \subseteq S \setminus\{ i \}}{\frac{v(L\cup \{ i \})-v(L)}{{{k} \choose |L|}}} +  \sum_{L \subseteq S \setminus\{ i \}}{\frac{v(L\cup \{ i,j \})-v(L\cup \{ j \})}{{{k} \choose |L|+1}}}  \right).
$$
Using that, by supermodularity, 
$$
v(L\cup \{ i \})-v(L) \leq v(L\cup \{ i,j \})-v(L\cup \{ j \}), 
$$ 
we get that 
$$
\phi_i(v_S)  \leq \phi_i(v_T),
$$
as claimed. Now, from this monotonicity, we conclude that
\beq \label{convcore}
v(S) = \sum_{i \in S}{\phi_i(v_S)} \leq  \sum_{i \in S}{\phi_i(v)},
\eeq
and therefore $\phi(v)$ is indeed in the core.
\epf

Combining Lemma \ref{supmod} and Lemma \ref{coreshap}, we immediately get:
\bc \label{coopcluethm} 
Let $v: 2^V \rightarrow [0,1] $ be a supermodular set function (cooperative game)  on the finite set $V$   with  $v(V)=1$, and let $\Gamma$ be a group acting on $V$ such that $v$ is invariant under the action of $\Gamma$. 
Then, for any $S \subseteq V$,
$$
v(S) \leq  \frac{|S|}{\min_{ i \in V}{|\Gamma \cdot i|}}.
$$
\ec

Recall that, with any function $f$ on a product space, we can associate a game $v_{f}$ by $ v_{f}(U): = \Var[\E[f \,|\, \mathcal{F}_U]]$. We get another  cooperative game if we define the information we gain about $Z=f(X_V)$ via information theoretic concepts (see Definition \ref{Infcluedef}):
$$
v^{I}_{f}(S)= I( Z : X_S ). 
$$
It is not difficult to see that, for product measures, both  $v_{f}$  and $v^{I}_{f}$ are convex games. The entropy version is immediate from the submodularity of entropy, which can be written as:
$$
    - H(X_{S}|Z) - H(X_{T}|Z) \leq  - H(X_{S \cap T}|Z)  - H(X_{S \cup T}|Z).
$$
Using that, for independent variables, the submodularity inequality is sharp, we get
\begin{multline*}
H(X_{S}) - H(X_{S}|Z) +  H(X_{T}) - H(X_{T}|Z) \leq \\
H(X_{S \cap T}) - H(X_{S \cap T}|Z) +  H(X_{S \cup T}) - H(X_{S \cup T}|Z).    
\end{multline*}
   
For the $L^2$ version, the supermodularity of $\Var(\E[f \,|\, \mathcal{F}_U])$ follows easily from the spectral description. Here we present an argument that does not require Fourier-Walsh expansion or Efron-Stein decomposition.

\bpr \label{clueconvex}
Let $f: \Omega^{V} \longrightarrow \R$, where $\Omega^{V}$ is endowed with a product measure. The set function (cooperative game) $v_f(S) = \Var(\E[f \,|\, \mathcal{F}_S])$ for ($S \subseteq V$) is supermodular (convex).
\epr

\bpf
First observe that, if $S \subseteq {\tilde{S}}$, then $\E[ \E[f \,|\, \mathcal{F}_{\tilde{S}}]\,|\, \mathcal{F}_S] = \E[f \,|\, \mathcal{F}_S]$ by the tower rule; then, using that conditional expectation is an orthogonal projection, we get that
$$
\Var(\E[f \,|\, \mathcal{F}_{\tilde{S}}]) - \Var(\E[f \,|\, \mathcal{F}_S]) = \Var(\E[f \,|\, \mathcal{F}_{\tilde{S}}] - \E[f \,|\, \mathcal{F}_S]).
$$
Therefore, the supermodularity condition \eqref{supermodular} can be rewritten as
\beq \label{l2supermod}
 \Var(\E[f \,|\, \mathcal{F}_T] -\E[f \,|\, \mathcal{F}_{S \cap T}]) \leq \Var(\E[f \,|\, \mathcal{F}_{S \cup T}] -\E[f \,|\, \mathcal{F}_{S}]).
\eeq

We will now use that the underlying measure is a product measure --- this will be the only place in the proof. Namely, for any $S \subseteq {\tilde{S}} \subseteq V$, using Lemma~\ref{keystep} for $({\tilde{S}} \setminus S)^{c}$ and ${\tilde{S}}$, we get
$$
\E[f \,|\, \mathcal{F}_S] = \E[\E[f \,|\, \mathcal{F}_{({\tilde{S}} \setminus S)^{c}}] \,|\, \mathcal{F}_{{\tilde{S}}}].
$$
This identity allows us to write 
\begin{align*}
\E[f \,|\, \mathcal{F}_{S \cap T}] &= \E[ \E[f \,|\, \mathcal{F}_{(T \setminus (S \cap T))^{c}}]\,|\, \mathcal{F}_T],\\
\E[f \,|\, \mathcal{F}_{S}] &= \E[ \E[f \,|\, \mathcal{F}_{((S \cup T) \setminus S)^{c}}]\,|\, \mathcal{F}_{S \cup T}].
\end{align*}
Since $T \setminus (S \cap T) =( S \cup T) \setminus S = T \setminus S$, \eqref{l2supermod} becomes
$$
 \Var(\E[f -\E[f \,|\, \mathcal{F}_{(T \setminus S)^{c}}]  \,|\, \mathcal{F}_T]) \leq \Var(\E[f -\E[f \,|\, \mathcal{F}_{(T \setminus S)^{c}}]  \,|\, \mathcal{F}_{S \cup T}]),
$$
which always holds, because orthogonal projection cannot increase the variance ($L^2$-norm).
\epf

We have shown that, for any real function $f$, the games $v_f$ (Proposition~\ref{clueconvex}) and  $v^{I}_f$ are convex whenever the underlying measure is a product measure; thus Corollary~\ref{coopcluethm} can be seen as a unified proof of Theorem~\ref{cluethm2} and (a generalization of) Theorem~\ref{entrcluethm}.  

Apart from product measures, however, we do not know any other example where the  condition of supermodularity of $v_f$ is satisfied. We note that in order to show that there is no sparse reconstruction for transitive functions  under some sequence of measures it would suffice to show the weaker condition that for some fixed $ k \in \N$  the set function $v^{k}_f(S) := (\clue(f\;|\;S))^{k}$  is supermodular. Indeed, in that case, by Corollary \ref{coopcluethm} we would get that $\clue(f\;|\;S)\leq (|U|/n)^{1/k}$.  

Observe that, for a transitive game with a non-empty core, the Shapley value, i.e., the uniform vector, will always be in the core. This is because the core is convex and itself is invariant under the group action. Therefore, one could weaken the condition of Proposition~\ref{coreshap} by only requiring the non-emptiness of the core. A classical result in Cooperative Game Theory (see, for example, \cite[Theorem 2.4]{BDT}) gives necessary and sufficient conditions for this. It has to be said, however, that on a practical level, the conditions of this theorem are not easy to verify.

\bth [Bondareva-Shapley]
The core of the game $v$ is non-empty if and only if  for every $\alpha: 2^V\setminus \emptyset \rightarrow [0,1]$ such that for every $i \in V$ 
$$
\sum_{S \subseteq V\;: \;i \in S}{\alpha(S)} = 1
$$
it holds that
$$
\sum_{S \subseteq V}{\alpha(S)v(S)} \leq v(V).
$$
\eth

\section{Sparse reconstruction for planar percolation}\label{s.perc}

\subsection{Background}\label{perctheory}

Bernoulli bond (or site) percolation at level $p$ on a graph $G$ means the random graph obtained by deleting every edge (or vertex) of a graph with probability $1-p$, independently. Here we only mention some basic concepts and results. For an introduction to  percolation theory,  criticality and other concepts, we advise the reader to consult \cite{Gr} in general, \cite{We} in two dimensions, and \cite{GS} with a focus on noise sensitivity questions.

In case $G$ is infinite,  we are interested  whether for a particular value of $p$ the arising random graph contains an infinite connected component. A simple coupling argument shows that this event is monotone increasing in $p$ and thus we introduce the critical value $p_c = \inf \left\{p : \mathbb{P}_p ( \exists \; \infty \; \text{cluster}) = 1 \right\}$.

Throughout this section we consider critical Bernoulli edge percolation on the square lattice (so every edge is open with probability $p=1/2$, independently). Our main focus will be the left-to-right crossing event  $\Cr_n$  on the $n\times (n-1)$ rectangle. This is the event that there exists  a path consisting of open edges between two vertices  located at opposite (left and right) sides of the rectangle. 

\begin{figure}[htbp]
\centerline{
\includegraphics[scale=0.25]{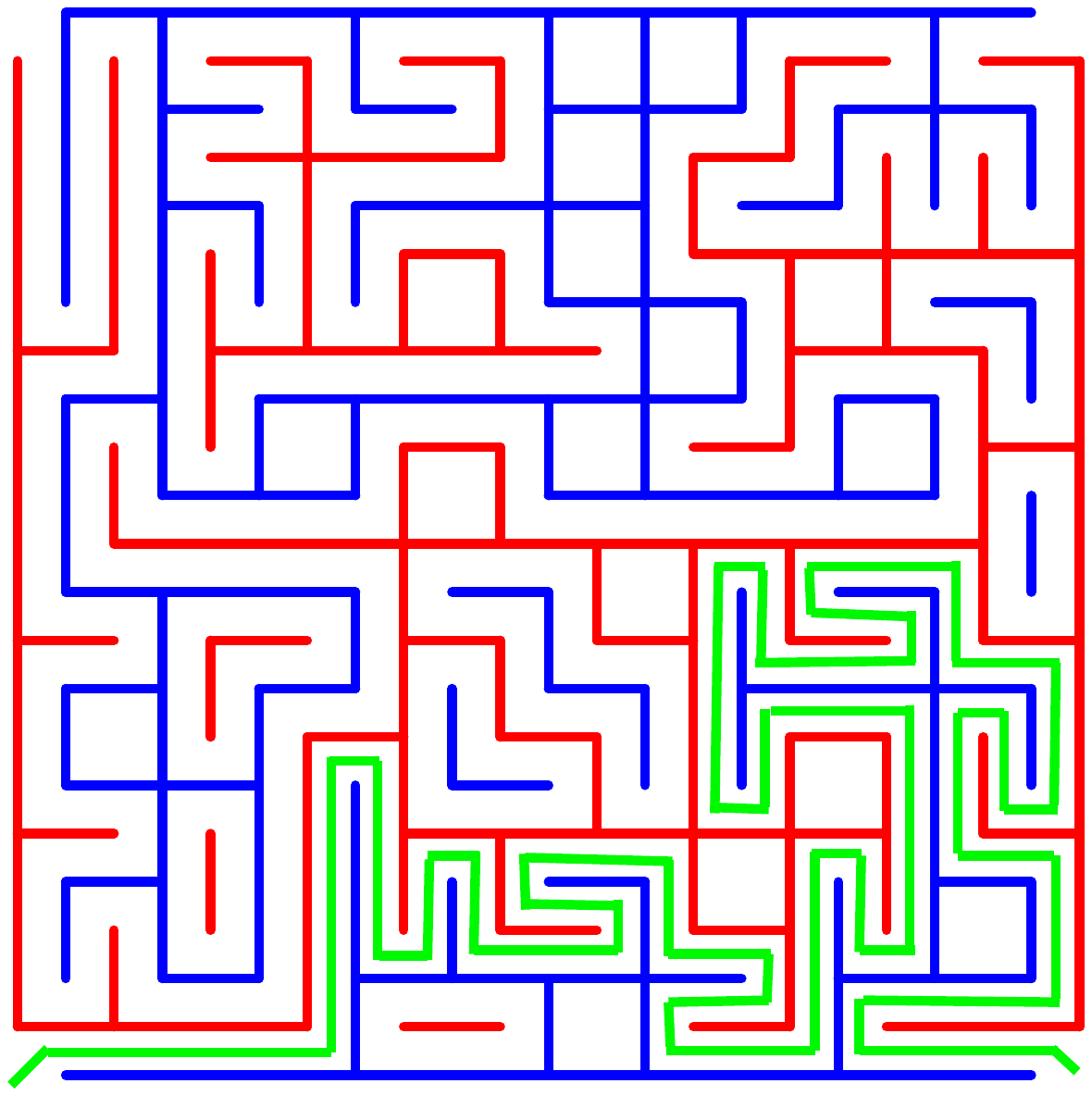}
}
\caption{A primal configuration (in red) where $\Cr_n$ is satisfied, its dual configuration (in blue), and the  exploration interface (in green) starting from the lower left corner.} 
\label{f.dual}
\end{figure}

It is known that if $p = \frac{1}{2}$, then $\p[\Cr_n]=\frac{1}{2}$ for any $n$. The reason behind this is that every percolation configuration on the square lattice induces a percolation configuration on the dual lattice: dual sites are the faces of the primal lattice, and two faces are connected in the dual configuration if the two faces are bordered by an edge which is closed in the primal percolation. The $n\times (n-1)$ rectangle $R_n$ has the important property that it is isomorphic to its dual graph, and there is a left-right crossing in the original lattice if{f} there is no top-bottom crossing of dual-open edges, and thus  $\p[\Cr_n]=\frac{1}{2}$ follows. To determine if $\Cr_n$ occurs or not, one can proclaim an open boundary condition on the left and right sides of $R_n$ and a dual-open boundary condition on the top and bottom sides, then start a so-called exploration interface from the lower left corner, going along the boundary of the open cluster attached to the left side of $R_n$, with open edges on the left of the exploration path, dual-open dual edges on the right. $\Cr_n$ occurs if{f} this exploration path hits the right side of $R_n$ and exits at the lower right corner, instead of hitting the top side and exiting at the upper left corner.  See Figure~\ref{f.dual}. 

The scale invariance $\p[\Cr_n]=\frac{1}{2}$ at $p=1/2$ is the crucial observation that suggests $p_c = \frac{1}{2}$, which is in fact known as the Harris-Kesten theorem; see the textbooks mentioned above.

It tuns out that the critical model in many graphs displays interesting, fractal-like features. There is a universality principle coming from statistical physics which connects the behaviour of various models around their phase transition. For example, physicists believe that percolation on any ``nice'' planar lattice $G$,  at the critical point $p_c(G)$, describes the same ``ideal'' percolation, only in possibly different frames.



We now introduce the so-called arm events some of which play an important role in the proof of our result. We start with the 1-arm event $A_1(R)$ on  $\Z^2$ --- we only consider this lattice, but the arm events can be defined in a similar way for any planar lattice. We think about the graph $\Z^2$ as embedded into the Euclidean plane $\R^2$ in the natural way together with its dual graph: we identify the vertices of the dual graph with the midpoints of the squares. Thus the dual can be identified with $\Z^2 + (1/2, 1/2)$. Let $B_{\mathbf{0}}(R)$ denote the (Euclidean) circle with center  $\mathbf{0}$ and radius $R$. Now for  $R>r>1$ we define the event  $A_1(r,R)$ as follows: There is a path of open edges between a vertex  $\mathbf{u} \in \Z^2 \cap B_{\mathbf{0}}(r)$  and some other vertex $\mathbf{v} \in \Z^2 \setminus B_{\mathbf{0}}(R)$. 

In a similar way other arm events may be defined. Our primary interest is the event $A^{+}_3(r, R)$, the three--arm event in a half plane. This is the event that there are three pairs of vertices  $(\mathbf{u}_i, \mathbf{v}_i), \; i  = 1,2,3$ in the upper half plane $H = \R \times \R^{+}$ in such a way that there are simultaneous paths of open edges $P_1$  connecting some vertex $\mathbf{u}_1 \in  B_{\mathbf{0}}(r) \cap H$ , with a vertex $\mathbf{v}_1 \in (B_{\mathbf{0}}(R))^{c}\cap H$ and $P_3$ connecting another vertex $\mathbf{u}_3 \in B_{\mathbf{0}}(r)\cap H$, with  some $\mathbf{v}_3 \in (B_{\mathbf{0}}(R))^{c}\cap H$, respectively. So far, these are only two arms: the third arm is a dual path between two dual vertices
 $\mathbf{u}_2 \in  B_{\mathbf{0}}(r)\cap H$, and  $\mathbf{v}_2 \in  (B_{\mathbf{0}}(R))^{c}\cap H$ that separates the two open paths $P_1$ and  $P_3$. Equivalently,  $P_1$ and  $P_3$ are not connected via open paths in the half-circle  $B_{\mathbf{0}}(R)\cap H$.
 
We shall also need $A^{+}_2(r, R)$, the two--arm event in the half plane. One  of the `arms'  is  a path of open edges connecting some vertex in the $r$-neighbourhood of $\mathbf{0}$ with another one outside the  $R$-neighbourhood of $\mathbf{0}$. The other `arm' is a dual one, which connects dual vertices in the respective neighbourhoods with dual open edges. Naturally, both of the paths are only allowed to use edges with both endpoints in $H$. Finally, we define $A^{++}_2(r, R)$, the two--arm event in the quarter plane  $Q := \R^{+} \times \R^{+}$. This is the same as $A^{+}_2(r, R)$, except that the two paths in question need to be contained in $Q$.
 
An important property of the critical model is that the probability of such arm events decay polynomially in $\frac{r}{R}$. Finding the exponent for the probability of some arm events are among the central questions of the field. The exponents for $A^{+}_2(r, R)$ and  $A^{+}_3(r, R)$, however, can be found with a combinatorial argument and in particular, they are known for the  $\Z^2$ lattice:

\bpr[{\cite[Appendix A]{LSW}}, {\cite[First Exercise Sheet]{We}}]\label{3armexponent}
For the $\Z^2$ lattice, 
$$
\alpha^{+}_2(r, R) := \p[A^{+}_2(r, R)]  \asymp \left( \frac{r}{R}\right), 
$$
and
$$
\alpha^{+}_3(r, R) := \p[A^{+}_3(r, R)]  \asymp \left( \frac{r}{R}\right)^{2}. 
$$
\epr

\subsection{No sparse reconstruction for critical planar percolation}\label{crnosr}

In this subsection we are going to prove Theorem~\ref{t.clueperc}, which we will now restate.

Any edge percolation configuration can be naturally identified with an  $\omega \in \{-1,1\}^{E}$, where $E$ is the edge set of the graph on which we percolate. In our case, $ \Cr_n :  \{-1,1\}^{E(R_n)}  \longrightarrow \{-1,1\}$ is the indicator function of the left-to-right crossing event in the $n \times (n-1)$ rectangle $R_n$.  We consider the critical probability $p=1/2$, thus we have the uniform measure on $\{-1,1\}^{E(R_n)}$ which we shall denote by $\p_n$.  Our result is the following:

\bth \label{clue_perc}
Let $ \Cr_n :  \left(\{-1,1\}^{E(R_n)}, \p_n \right) \longrightarrow \{-1,1\}$ be the left-right crossing event as above.  Then there is no sparse reconstruction for $\Cr_n$, that is for any sequence $U_n \subset E(R_n)$ with $ \lim_{n \to \infty}{|U_n|/|E(R_n)|} =0$, we have
$$
 \lim_{n \to \infty}{\clue(\Cr_n\;|\;U_n)} = 0.
$$
\eth

Here is a brief summary of what we are going to do. First of all, we embed $R_n$ into the torus graph $T_n=\Z_n \times \Z_{n-1}$, which is just the lattice $\Z^2$ quotiented by the natural translation action of the subgroup $n\Z \times (n-1)\Z$; this quotient has the same vertex set as $R_n$, and some extra edges ``across the boundary''. The group $\Z_n \times \Z_{n-1}$ acts transitively on the vertices of $T_n$ by translations, while the action on the edges has two orbits: the horizontal and the vertical edges. In particular, the function  $\Cr_n$  has natural translations $\Cr_n^\mathbf{t}:  \{-1,1\}^{E(T_n)}  \longrightarrow \{-1,1\}$ for any $\mathbf{t} \in \Z_n \times \Z_{n-1}$.

The key percolation ingredient will be Lemma~\ref{continuity}, which says that there is an absolute constant $K$ such that, for every $\delta>0$ and large enough $n \in \N$, if $\mathbf{t} \in \Z_n \times \Z_{n-1}$ has length at most $\delta n$, then the correlation between $\Cr_n$ and $\Cr^{\mathbf{t}}_n$ is at least $1-K \delta$. Now Lemma~\ref{proj} tells us that if two functions are highly correlated and one of them has high clue with respect to a subset, then the other one also has high clue with respect to the same subset. Thus, if we assume that there is a sparse sequence of subsets $\{U_n\}_{n \in \N}$ such that $\clue(\Cr_n\;|\;U_n)>c$ for infinitely many $n \in \N$, then we also have $\clue(\Cr_n\;|\;U_n^{\mathbf{t}})=\clue(\Cr^{-\mathbf{t}}_n\;|\;U_n)>c/2$ for all $\|\mathbf{t}\|_\infty \leq \delta n $, for some $\delta>0$ chosen according to $c$. 



From here, we will follow Proposition~\ref{l.almosttran}. Let us define the sublattice
\begin{equation}\label{e.Ldelta}
L _{\delta} : = \big\{\pm\floor{\delta n}, \pm3\floor{\delta n}, \dots, \pm (2N-1)\floor{\delta n}\big\}^2 \subset \Z_n \times \Z_{n-1}\,,
\end{equation}
where $N$ is the largest integer with $(4N-2)\floor{\delta n} \leq n-1$. It is clear that $L _{\delta}$  has on the order of $1/\delta^2$ elements, which is a constant independent of $n$. Therefore, if we define $W_{\delta,n}$ to be the union of all the $L _{\delta}$-translates of $U_n$, that is still sparse. Moreover, every translation $\gamma \in \Z_n \times \Z_{n-1}$ can be written as a sum of an $\ell \in L_\delta$ and some vector $\|\mathbf{d}\|_\infty \leq \delta n$, which implies that $W_{\delta,n}^\gamma$ contains $U^\mathbf{d}$, and hence $\clue(\Cr\;|\;W_{\delta,n}^\gamma) > c/2$. Taking a uniformly random $\gamma \in \Z_n \times \Z_{n-1}$, the resulting random set $\cW:=W_{\delta,n}^\gamma$ will have small revealment, but expected clue at least $c/2$, contradicting Theorem~\ref{t.cluerandom}.

We now do this in more detail. Let $0<\delta<1$. For a $\mathbf{t} \in \Z_n \times \Z_{n-1}$, we will denote the rectangle $\mathbf{t} + [-\floor{\delta n}, \floor{\delta n}] ^2\subset \Z_n \times \Z_{n-1}$ by $R_\delta(\mathbf{t})$. It is straightforward to see that $(2\delta n - 1)^2\leq |R_\delta(t)| \leq (2\delta n + 1)^2$.  

The following lemma might already be somewhere in the percolation literature; for instance, \cite[Appendix A]{SSG} gives a more general result, with a weaker bound.

\bl  \label{continuity}
Let $R_\delta:= R_\delta(\mathbf{0}) =  [-\floor{\delta n},\floor{\delta n}] ^2 $ as above. Then there is an absolute constant $K>0$ such that, for every $\mathbf{d} \in R_\delta$, 
$$ 
\Corr(\Cr_n , \Cr_n ^{\mathbf{d}}) \geq 1- K\delta,
$$
for every $n$ large enough.
\el

\bpf
Let $\mathbf{d} \in R_\delta$. We will show that 
$$\p [\Cr_n \neq  \Cr_n ^{\mathbf{d}}] \leq O (\delta);$$
from this, the statement of the lemma follows, since
$$
 \Corr(\Cr_n , \Cr_n ^{\mathbf{d}})= 1- 2\p [\Cr_n  \neq  \Cr_n ^{\mathbf{d}}]. 
$$

\begin{figure}[htbp]
\vskip 0.3 cm
\SetLabels
(-0.05*-0.06){$(0,0)$}\\
(0.85*-0.06){$(n-\delta n,0)$}\\
(1.05*1.03){$(n,n-1)$}\\
(1.02*0.1)$\Bigg\} r$\\
(1.04*0.25) $B$\\
(1.04*0.45) $A$\\
\endSetLabels
\centerline{
\AffixLabels{
\includegraphics[scale=0.5]{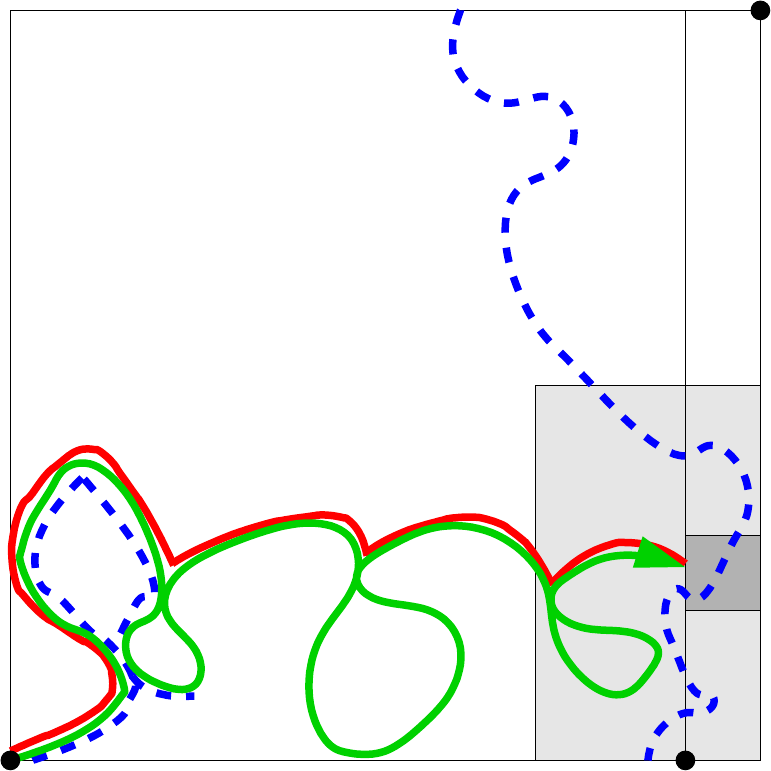}
}}
\vskip 0.4 cm
\caption{If $\Cr_n \neq \Cr_n ^{(-\delta n,0)}$, then the primal (red) left-to-right crossing and the dual (dashed blue) top-to-bottom crossing create a 3-arm event in a half-plane from $\partial B$ to $\partial A$, and a 2-arm event in a quarter-plane from $\partial A$ to $\partial R_n$, for the $\delta n \times \delta n$ box $B$ that is first hit by the exploration interface.} 
\label{f.3arm}
\end{figure}

As an extreme case, let us first assume that $\mathbf{d} =(-{\delta n},0)$, ignoring the integer part from $\floor{\delta n}$ for easier notation. We are going to show that the event $\left\{\Cr_n \neq  \Cr_n ^{\mathbf{d}} \right\}$ basically entails a 3-arm event in a half plane from the boundary of one of $O(1/\delta)$ possible $\delta n \times \delta n$ boxes to a distance of order $n$, which has, by Proposition~\ref{3armexponent}, a probability $O(1/\delta) \, O(\delta^2) = O(\delta)$ to occur. See Figure~\ref{f.3arm} to have a rough idea before the detailed proof.

By symmetry, we may assume that $\Cr_n$ does not hold, while $\Cr_n ^{(-\delta n,0)}$ does. The latter event implies that there is a left-to-right crossing in the smaller rectangle with lower left corner $(0,0)$ and upper right corner $(n-\delta n,n-1)$. Start an exploration interface from $(0,0)$ with open edges on the left and dual-open dual edges on the right, as in Figure~\ref{f.dual}, until it first hits the vertical line $(n-\delta n,\cdot)$; at this moment, along the left boundary of the interface we have found the lowest open left-to-right crossing from $(0,\cdot)$ to $(n-\delta n,\cdot)$. 

Let us divide the rectangle with opposite corners $(n-\delta n,0)$ and $(n,n-1)$ into boxes of size close to $\delta n \times \delta n$, and let $B$ be the box that is hit by the exploration interface (darker grey on the figure). Since we have assumed that $\Cr_n$ does not hold, there must be a top-to-bottom dual crossing of $R_n$, and because of the open crossing we have already found, this dual crossing must go through $B$. Now let the distance of $B$ from the top and bottom boundaries of $R_n$ be $r$, and let $A$ be the half-square of radius $r+\delta n/2$ centered at $B$ (lighter grey on the figure). Then the primal and dual crossings create a half-plane 3-arm event from distance $n\delta$ to $r$, and a quarter-plane $2$-arm event from a distance of order $r$ to $n$. For the quarter-plane probabilities $\alpha^{++}_{2}(r,R) : = \p[A^{++}_{2}(r,R)]$ we can use the obvious bound 
$$
\alpha^{++}_{2}(r,R) \leq \alpha^{+}_{2}(r,R).
$$
Using a dyadic division $r\in [n\delta 2^j, n\delta 2^{j+1})$ for the possible values of $r$, and using that at each scale there are order $2^j$ possible locations for $B$, we have, by the arm probability bounds in Proposition~\ref{3armexponent}:
\begin{align*}
 \p [\Cr_n \neq  \Cr_n ^{\mathbf{d}}] &\leq O(1) \sum_{j=1}^{\log (1/\delta)} \alpha^{+}_{3}(n\delta,n\delta 2^j) \, \alpha^{++}_{2}(n\delta 2^j,n)\, 2^j \\
 &\leq O(1) \sum_{j=1}^{\log (1/\delta)} (2^j)^{-2} \, (\delta 2^j)^{2} 2^j\\
 & \asymp \delta^2  \sum_{j=1}^{\log (1/\delta)} 2^j \asymp \delta\,,
\end{align*}
as claimed.

When $\mathbf{d} =(t,0)$ for some $-\delta n \leq t \leq \delta n$, we similarly have the upper bound $O(t)=O(\delta)$. 

In case $\mathbf{d} =(0,t)$ we have exactly the same argument, exploiting duality and the $\pi/2$ rotational symmetry of the model: since $\Cr_n$ does not happen i{f}f there is a dual up-down crossing, the event $\{\Cr_n \neq  \Cr_n^{\mathbf{d}}\}$ means that the two vertically $O(\delta)$-shifted rectangles disagree about having dual up-down crossings or not, which is exactly as unlikely as horizontally $O(\delta)$-shifted rectangles disagreeing about left-to-right primal crossings. 

The case of a general $\mathbf{d} \in  R_\delta$ now easily follows. If  $\left\{\Cr_n  \neq  \Cr_n ^{\mathbf{d}} \right\}$, then  either $\big\{\Cr_n  \neq  \Cr_n ^{\mathbf{d}_x} \big\}$  or $\big\{\Cr_n \neq  \Cr_n ^{\mathbf{d}_{y}} \big\}$, where $\mathbf{d}_x$ and $\mathbf{d}_{y}$ are the projections of $\mathbf{d}$ onto the first and the second coordinates, respectively. 

As a consequence, $\p [\Cr_n  \neq  \Cr_n ^{\mathbf{d}}]\leq \p [\Cr_n  \neq  \Cr_n ^{\mathbf{d}_x}] + \p [\Cr_n  \neq  \Cr_n ^{\mathbf{d}_{y}}] \leq O (\delta) $.
\epf

\bpf[Proof of Theorem \ref{clue_perc}]
Let $U_n \subseteq \Z_{n}^2 $ be a sparse sequence of subsets; i.e.,  $\lim_{n} {\frac{|U_n|}{n^2}} =0$. 

Recall that for any $\delta>0$ we have the notation $R_\delta=  [-\floor{\delta n},\floor{\delta n}] ^2$, and the sublattice $L _{\delta}$ from~(\ref{e.Ldelta}). Clearly, $N\leq 1/(2\delta)$ if $\delta$ is small and $n$ is large enough, hence $|L_\delta|=4 N^2 \leq 1/\delta^2$.

It is straightforward to check that any $\mathbf{t} \in \Z_n \times \Z_{n-1}$ can be written as the sum of a $\mathbf{d} \in R_\delta$ and an $\mathbf{\ell} \in L_\delta$. Recall that, according to Lemma~\ref{continuity}, 
$$ 
\Corr(\Cr_n , \Cr_n ^{\mathbf{d}}) \geq 1- K\delta,
$$
for any  $\mathbf{d} \in R_\delta$. Therefore, by Proposition~\ref{l.almosttran}, using that the action of $\Z_n \times \Z_{n-1}$ on $ E(\Z_n \times \Z_{n-1})$ has two orbits, the set of vertical edges $E_v$ and the set of horizontal edges $E_h$, 
$$ 
\clue(\Cr_n\;|\;U_n)\leq  \frac{\left|L _{\delta}\right|\left|U_n\right|}{n(n-1)} +5\sqrt{K \delta} = \frac{1}{\delta^2} \frac{\left|U_n\right|}{n(n-1)} + 5\sqrt{K \delta}.
$$ 
By choosing $\delta$ small, and then $n$ large enough, and using our assumption $\left|U_n\right|/n^2 \to 0$, both terms in the last upper bound can be made arbitrarily small, hence the theorem follows.
\epf

\section{Some open problems}\label{s.open}

Here we collect the open problems raised somewhere in the paper and some further ones. We have not thought thoroughly about these questions, but we would definitely be interested in the answers.

From Subsection~\ref{siginf}, we have the following question relating significance and influence; see the discussion there.

\begin{quest}\label{q.siginf}
Characterise sequences of Boolean functions  such that for any sequence of subsets with constant density $I(f_n \;|\;U_n) \gg \sig(f_n \;|\;U_n)$ holds, or where $I(f_n \;|\;U_n) \asymp \sig(f_n \;|\;U_n)$, respectively. 
\end{quest}

From Section \ref{s.infoclue}, comparing $L^2$-clue and I-clue:

\begin{quest}\label{q.L2I}
Is there a sequence  of functions $f_n \in L^2( \{-1,1\}^{V_n},\pi_n^{\otimes V_n})$  and a corresponding sequence of subsets $U_n \subseteq V_n$ such that 
\begin{enumerate}
    \item $\clue^{I}(f_n\;|\;U_n) \ll \clue(f_n\;|\;U_n)$
    \item $\clue^{I}(f_n\;|\;U_n) \gg \clue(f_n\;|\;U_n)$?
\end{enumerate}
What is the answer if we ask $f_n$ to be Boolean for all $n \in \N$?
What is the answer if we allow for non-product measures on $\{-1,1\}^{V_n}$?
\end{quest}


\def\cp{\mathsf{cp}}
\def\rcp{\mathsf{rcp}}

Thinking primarily of the setup of iid measures $\pi^{\otimes V}$, let us define the clue profile:
$$
\cp_f(\delta) := \sup\big\{   \clue(f \given U)  : |U|/|V| \leq \delta \big\}.
$$
We would like to define a random clue profile, as well, for sparse random sets $\mathcal{U}$ that are independent of the $\sigma$-algebra of $\pi^{\otimes V}$. However, as mentioned in the discussion before Theorem~\ref{NoRSR}, the simple strategy of asking all bits with probability $\delta$ achieves the maximal average clue among all possible random sets of revealment $\delta$. To come around this, we only allow random subsets of a fixed size:
\begin{align*}
\rcp_f(\delta) := \sup \big\{ \E\big[\clue(f \given \cU)\big] : \cU \textrm{ is a random subset with } |\cU|=\lfloor\delta n\rfloor \textrm{ and revealment}\leq\delta \}.
\end{align*}
\begin{quest}\label{q.cpshape} 
Can we say anything general about the shapes of these monotone increasing functions $\cp_f: [0,1] \lora [0,1]$ and  $\rcp_f: [0,1] \lora [0,1]$, beyond $cp_f(\delta)\leq \delta$ that comes from Theorem~\ref{t.cluegen} and Theorem~\ref{NoRSR}? Almost anything is possible, as in \cite{ASP}? When do we have a sharp threshold, as for monotone graph properties \cite{hunt} or in the cutoff phenomenon for random walk mixing times \cite{AD,BHP}?
\end{quest}


In particular, one might be interested in the worst possible cases.

\begin{quest}\label{q.rcplow} 
For a balanced sequence of monotone functions $f_n \in L^2( \{-1,1\}^{V_n},\pi_n^{\otimes V_n})$ how small $\rcp_{f_n}(\delta)$ can be for some fixed $\delta \in (0,1)$? In particular, is it possible that for some $\delta$ we have  $\lim_{n}{\rcp_{f_n}(\delta)} = 0$?
\end{quest}

Monotonicity is necessary because of the parity $\chi_{V_n}$.  

As one example, for $\Cr_n$ on the triangular lattice, one can try a (randomly positioned) subsquare of size $\lfloor\sqrt{\delta} n\rfloor \times \lfloor \sqrt{\delta} n \rfloor$. Such a subset, as calculated in \cite{GPS}, has $\clue$ of order $\delta^{5/4+o(1)}$, and one might expect that this is the best one can achieve. 

Looking at the proofs of Theorems~\ref{cluethm} and~\ref{NoRSR}, one may speculate that for noise-sensitive sequences, where the spectral sample $\Spec_n$ is typically much larger than just one element, these results should be very far being sharp, and therefore $\clue(f_n \given \cU_n) \ll \delta(\cU_n)$ should always hold. However, the noise-sensitive sequence of the $\Tr$ functions is a counterexample: if we choose a uniformly random $\delta$ proportion of the tribes, and ask all the bits in those, then the clue of this subset will be of order $\delta$.

The Schramm-Steif small revealment noise sensitivity theorem \cite{SS} does not have a converse: there are monotone noise-sensitive functions without a small revealment algorithm \cite[Section 8.6]{GS}. The following question would aim at a converse in terms of clue. Recall that $\B^\delta$ denotes the random set that contains every  $v \in V$ independently, with probability $\delta$.

\begin{quest}\label{q.cpsensemon} Is it true for every noise-sensitive  sequence of monotone Boolean functions $f_n$ that, for any fixed $\delta>0$, $\rcp_{f_n} (\delta) \gg \E\big[\clue(f_n \given \B^\delta)\big]$?
\end{quest}

We can extend Question \ref{q.cpsensemon} to general noise-sensitive Boolean sequences and more general sequences $\delta_n$: 

\begin{quest}\label{q.cpsense} Is it true for every noise-sensitive  sequence of  Boolean functions $f_n$ that, for any positive sequence $\{\delta_n\}_{n \in \N}$ satisfying 
$\frac{|\Spec_f|}{\delta_n n} \xrightarrow{\p}  0,$ 
we have $\rcp_{f_n} (\delta_n) \gg \E\big[\clue(f_n \given \B^{\delta_n})\big]$? 
\end{quest}

This question, on an intuitive level, is in the direction of the Fourier Entropy vs Influence conjecture of Friedgut and Kalai \cite{thresh} (see also \cite{Kalai}, \cite{ODWZ}, \cite{KKLMS}, and the references therein), which can be formulated by saying that the spectral sample of a Boolean function always has some structure that reduces its entropy, compared to a uniformly random set of a similar size. Question~\ref{q.cpsense}  also hints at such a structure: it should be possible to cover $\Spec_f$ more efficiently by some $\cU$ than by an i.i.d.~random set.

A final remark is that it could be interesting to find an exact connection between our Theorem \ref{t.cluegen} and \cite[Theorem 7]{log}.

\vfill\eject

\ \\
{\bf P\'al Galicza}\\
Alfr\'ed R\'enyi Institute of Mathematics\\
Re\'altanoda u. 13-15, Budapest 1053 Hungary\\
\texttt{galicza[at]gmail.com}\\
\ \\
{\bf G\'abor Pete}\\
Alfr\'ed R\'enyi Institute of Mathematics\\
Re\'altanoda u. 13-15, Budapest 1053 Hungary\\
and\\
Department of Stochastics, Institute of Mathematics\\
Budapest University of Technology and Economics\\
M\H{u}egyetem rkp.~3., Budapest 1111 Hungary\\
\texttt{robagetep[at]gmail.com}\\
\url{http://www.math.bme.hu/~gabor}\\

\end{document}